\documentclass[11pt, reqno]{article}
\usepackage[T1]{fontenc}
\usepackage{amsmath,amssymb,amsfonts,amsthm}
\usepackage{color}

\usepackage[colorlinks=true,linkcolor=blue,citecolor=blue,urlcolor=blue,pdfborder={0 0 0}]{hyperref}
\usepackage{cleveref}

\usepackage{caption}
\usepackage{thmtools}

\usepackage{mathrsfs}

\crefname{theorem}{Theorem}{Theorems}
\crefname{thm}{Theorem}{Theorems}
\crefname{lemma}{Lemma}{Lemmas}
\crefname{lem}{Lemma}{Lemmas}
\crefname{remark}{Remark}{Remarks}
\crefname{prop}{Proposition}{Propositions}
\crefname{defn}{Definition}{Definitions}
\crefname{corollary}{Corollary}{Corollaries}
\crefname{conjecture}{Conjecture}{Conjectures}
\crefname{question}{Question}{Questions}
\crefname{chapter}{Chapter}{Chapters}
\crefname{section}{Section}{Sections}
\crefname{figure}{Figure}{Figures}
\crefname{example}{Example}{Examples}

\theoremstyle{plain}
\newtheorem{thm}{Theorem}[section]

\newtheorem{lemma}[thm]{Lemma}
\newtheorem{theorem}[thm]{Theorem}

\newtheorem{corollary}[thm]{Corollary}

\newtheorem{prop}[thm]{Proposition}

\theoremstyle{definition}

\theoremstyle{remark}
\newtheorem{remark}[thm]{Remark}

\numberwithin{equation}{section}

\renewcommand{\P}{\mathbb P}
\newcommand{\E}{\mathbb E}

\newcommand{\R}{\mathbb R}
\newcommand{\Z}{\mathbb Z}
\newcommand{\N}{\mathbb N}

\newcommand{\cE}{\mathcal E}
\newcommand{\cF}{\mathcal F}

\newcommand{\cP}{\mathcal P}

\newcommand{\cW}{\mathcal W}

\usepackage[margin=1in]{geometry}
\usepackage{tikz}
\usepackage{pgfplots}
\usepackage{extarrows}
\usepackage{titling}
\usepackage{commath}
\usepackage{wasysym}
\usepackage{mathtools}
\usepackage{titlesec}
\newcommand{\eps}{\varepsilon}
\usepackage{bbm}
\usepackage{setspace}
\usepackage{enumitem}
\usepackage{booktabs}

\usepackage[textsize=tiny]{todonotes}

\usepackage{cite}

\newcommand{\Aut}{\operatorname{Aut}}

\newcommand{\bP}{\mathbf P}
\newcommand{\bE}{\mathbf E}

\newcommand{\Lip}{\operatorname{Lip}}

\def\P{\mathbb{P}}

\DeclareMathOperator*{\argmin}{argmin}

\DeclareMathOperator{\Gr}{Gr}

\newcommand{ \TV}{ \mathrm{TV} }

\DeclareMathSymbol{\leslant}{\mathalpha}{AMSa}{"36} 
\DeclareMathSymbol{\geslant}{\mathalpha}{AMSa}{"3E} 
\DeclareMathSymbol{\eset}{\mathalpha}{AMSb}{"3F}     

\renewcommand{\epsilon}{\varepsilon}

\newcommand{\defeq}{\vcentcolon=}
\newcommand{\eqdef}{=\vcentcolon}

\makeatletter
\pgfdeclareshape{slash underlined}
{
  \inheritsavedanchors[from=rectangle] 
  \inheritanchorborder[from=rectangle]
  \inheritanchor[from=rectangle]{north}
  \inheritanchor[from=rectangle]{north west}
  \inheritanchor[from=rectangle]{north east}
  \inheritanchor[from=rectangle]{center}
  \inheritanchor[from=rectangle]{west}
  \inheritanchor[from=rectangle]{east}
  \inheritanchor[from=rectangle]{mid}
  \inheritanchor[from=rectangle]{mid west}
  \inheritanchor[from=rectangle]{mid east}
  \inheritanchor[from=rectangle]{base}
  \inheritanchor[from=rectangle]{base west}
  \inheritanchor[from=rectangle]{base east}
  \inheritanchor[from=rectangle]{south}
  \inheritanchor[from=rectangle]{south west}
  \inheritanchor[from=rectangle]{south east}
  \inheritanchorborder[from=rectangle]
  \foregroundpath{
    \southwest \pgf@xa=\pgf@x \pgf@ya=\pgf@y
    \northeast \pgf@xb=\pgf@x \pgf@yb=\pgf@y
    \pgf@xc=\pgf@xa
    \advance\pgf@xc by .5\pgf@xb
    \pgf@yc=\pgf@ya
    \advance\pgf@xc by -1.3pt
    \advance\pgf@yc by -1.8pt
    \pgfpathmoveto{\pgfqpoint{\pgf@xc}{\pgf@yc}}
    \advance\pgf@xc by  2.6pt
    \advance\pgf@yc by  3.6pt
    \pgfpathlineto{\pgfqpoint{\pgf@xc}{\pgf@yc}}
    \pgfpathmoveto{\pgfqpoint{\pgf@xa}{\pgf@ya}}
    \pgfpathlineto{\pgfqpoint{\pgf@xb}{\pgf@ya}}
 }
}
\tikzset{nomorepostaction/.code=\let\tikz@postactions\pgfutil@empty}

\makeatother

\let\oldchi=\chi
\renewcommand{\chi}{\text{\raisebox{\depth}{\(\oldchi\)}}}

\newcommand{\diam}{\operatorname{diam}}

\pretitle{\begin{flushleft}\Large}
\posttitle{\par\end{flushleft}\vskip 0.5em
}
\preauthor{\begin{flushleft}}
\postauthor{
\\
\vspace{0.5em}
\footnotesize{
$\spadesuit$ Department of Mathematics, ETH Zurich
\\
Email: \href{mailto:camillo.brena@math.ethz.ch}{camillo.brena@math.ethz.ch}\\
$\heartsuit$ Department of Mathematics, Princeton University and
PMA, Caltech
\\
Email: \href{mailto:t.hutchcroft@princeton.edu}{t.hutchcroft@princeton.edu}, \href{mailto:t.hutchcroft@caltech.edu}{t.hutchcroft@caltech.edu}\\
$\diamondsuit$ 
Institute for Mathematics, Leipzig University, Germany
\\
Email: \href{mailto:cfmuench@gmail.com}{cfmuench@gmail.com}}
\end{flushleft}}
\predate{\begin{flushleft}}
\postdate{\par\end{flushleft}}
\title{{\bf Transitive graphs of non-negative Ollivier--Ricci curvature have polynomial growth}}

\renewenvironment{abstract}
 {\par\noindent\textbf{\abstractname.}\ \ignorespaces}
 {\par\medskip}

\titleformat{\section}
  {\normalfont\large \bf}{\thesection}{0.4em}{}

\author{{\bf Camillo Brena$^\spadesuit$, Tom Hutchcroft$^\heartsuit$, and Florentin M\"unch$^\diamondsuit$}}

\begin{document}

\maketitle

\abstract{We prove that every transitive graph of non-negative Ollivier--Ricci curvature has polynomial growth. We deduce from this and a theorem of Brena and Bru\`e that a finitely generated group admits a Cayley graph of non-negative Ollivier--Ricci curvature if and only if it is virtually abelian. Beyond the transitive setting, we also prove a quasi-polynomial $r^{O(\sqrt{\log r})}$ volume bound for balls around typical vertices in finite bounded-degree graphs of non-negative Ollivier--Ricci curvature, greatly strengthening a theorem of Salez (GAFA 2022).
}


\setstretch{1.1}

\section{Introduction}
This paper concerns the geometric and spectral properties of bounded-degree graphs with \emph{non-negative Ollivier--Ricci curvature}, a discrete analogue of non-negative Ricci curvature for Riemannian manifolds introduced in an inspiring series of works by Ollivier \cite{ollivier2007ricci,ollivier2009ricci,ollivier2010survey} and recently given renewed prominence in works such as \cite{munch2023mixing,pedrotti2025new,caputo2025entropy,caputo2024entropy,salez2025cutoff,salez2024spectral,salez2023cutoff,pedrotti2025transport,munch2024perpetual}; see \cite{salez2025modern} for a recent survey on applications to Markov chain mixing times and the cut-off phenomenon.




Let us now give the relevant definitions.
Recall that if $G=(V,E)$ is a locally finite graph without isolated vertices,
the \textbf{lazy random walk} on $G$ is the process $(X_n)_{n\geq 0}$ which at each step either stays in place with probability $1/2$ or else crosses a uniform random edge emanating from its current location, independently of everything it has done previously. 
 The \textbf{Markov operator} (a.k.a.\ \textbf{transition matrix}) of lazy random walk on $G$ is defined by setting 
\[P(x,y)=\mathbf{P}_x(X_1=y)=\frac{1}{2}\mathbbm{1}(x=y)+ \frac{\#\{e\in E^\rightarrow : e^-=x,e^+=y\}}{2\deg(x)},\] where $\mathbf{P}_x=\mathbf{P}_x^G$ denotes the law of the lazy random walk started at $x$ and where we think of each edge of $G$ as corresponding to a pair of oriented edges, one in each direction, with the set of oriented edges denoted by $E^\rightarrow$. (We allow our graphs to have self-loops and multiple edges.)

A locally finite graph $G=(V,E)$ is said to have \textbf{non-negative Ollivier--Ricci curvature} if it is non-empty, has no isolated vertices\footnote{These conditions are not standard but save us from the repeated consideration of trivialities.}, and any of the following equivalent conditions hold: 
\begin{enumerate}
  \item The Markov operator $P$ defines a contraction on the space of Lipschitz functions on $G$: If $f:V\to \R$ is any function then $\operatorname{Lip}(Pf) \leq \operatorname{Lip}(f)$ where $Pf(x)=\sum_{y\in V} P(x,y)f(y)$ and $\operatorname{Lip}(f)=\sup\{|f(x)-f(y)|:x\sim y\}$.
  \item 
The Markov operator $P$ defines a contraction on the ($L^1$) Wasserstein space $\cW$ associated to the graph $G$: If $\mu$ and $\nu$ are probability measures on $V$ then 
$W(\mu P,\nu P) \leq W(\mu,\nu)$,
 where we recall that if $\mu$ is a measure on $V$ then $\mu P$ is the measure defined by $\mu P(x) = \sum_{y\in V}\mu(y)P(y,x)$ and the Wasserstein distance between two measures $\mu$ and $\nu$ is defined by 
  \begin{align*}
   W(\mu,\nu) &=  \inf \Bigl\{ \E [d(X,Y)] : (X,Y) \text{ is a coupling of }\mu\text{ and }\nu \Bigr\}\\
   & =\sup\Bigl\{ \bigl|\sum_x f(x)\mu(x)-\sum_x f(x)\nu(x)\bigr| : f:V\to \R \text{ finitely supported},\, \operatorname{Lip}(f)\leq 1 \Bigr\};
 \end{align*}
 the equivalence between these two quantities is an instance of Kantorovich duality.
 \item The Wasserstein contraction property $W(\delta_x P,\delta_y P) \leq 1$ holds for all pairs of neighbouring vertices $x,y$, where $\delta_x$ denotes the Dirac mass at $x$.
\item The random walks $X$ and $Y$ started at two vertices $x$ and $y$ can be coupled in such a way that $(d(X_n,Y_n))_{n=0}^\infty$ is a supermartingale with respect to its natural filtration, meaning explicitly that
\[\E \Bigl[d(X_{n+1},Y_{n+1}) \;\Big|\; (X_i)_{i=0}^n,(Y_i)_{i=0}^n \Bigr]\leq d(X_n,Y_n)\]
almost surely for every $n\geq 0$, where $\E$ denotes an expectation taken with respect to the coupled pair of random walks. 
\end{enumerate}
See \cite[Theorem 4.1]{salez2025modern} for a proof that these conditions are equivalent. 
This definition is motivated by a theorem of von Renesse and Sturm \cite{von2005transport} characterising Riemannian manifolds of non-negative Ricci curvature as precisely those for which the heat semigroup defines a contraction on the space of Lipschitz functions. Examples of non-negatively curved graphs include Cayley graphs of abelian groups and, more generally, Cayley graphs of any group taken with respect to a conjugation-invariant generating set \cite[Theorem 4.6]{salez2025modern}.
 Non-negative (and positive) Ollivier--Ricci curvature also occurs in many examples of \emph{interacting particle systems}, where it is known as \emph{the Dobrushin criterion} \cite{dobrushin1985constructive} and has been studied since the 1980s.

Just as non-negative curvature is a local (i.e.\ ``infinitesimal'') property of Riemannian manifolds, non-negative Ollivier--Ricci curvature is a local
 property of graphs in the sense that it is determined by the isomorphism classes of balls of radius $2$ appearing in the graph.
We are interested in the 
\emph{large-scale} geometric consequences of this local property. This is a classical problem in the Riemannian context, where there are many important and now-classical theorems concerning the large-scale geometry of manifolds of non-negative Ricci curvature. One of the most significant basic theorems in this direction is the \emph{Bishop--Gromov inequality} \cite{petersen2006riemannian}, which states that every Riemannian $n$-manifold  $M$ of non-negative Ricci curvature has volume growth bounded by that of $\R^n$, and in fact that
$r^{-n}\operatorname{Vol}(B(x,r))$
is monotone decreasing in $r$ for every $x\in M$.

In his influential survey \cite{ollivier2010survey}, Ollivier \cite[Problem L]{ollivier2010survey} asked whether there is any analogue of the Bishop--Gromov inequality for graphs of non-negative Ollivier--Ricci curvature. Progress in this direction has so far been extremely limited. Indeed, it was proven only recently in a breakthrough work of Salez \cite{salez2022sparse} that bounded-degree graphs of non-negative Ollivier--Ricci curvature \emph{are not expanders}, a qualitative property vastly weaker than the polynomial growth ensured by the Bishop--Gromov inequality. Subsequent work of the second named author and Lopez \cite{hutchcroft2024relation} established a quantitative version of Salez's theorem, but this quantitative version remained too weak even to conclude that bounded-degree graphs of non-negative Ollivier--Ricci curvature have subexponential volume growth. Moreover, while Salez proved his result by establishing that the random walk on a bounded-degree graph of non-negative Ollivier--Ricci curvature satisfies a sublinear expected displacement bound of the form $\mathbf{E}_{X_0\sim \pi} [d(X_0,X_t)] =o(t)$, both his methods and those of \cite{hutchcroft2024relation} did not yield any quantitative version of this estimate. Since Brownian motion on manifolds of non-negative Ricci curvature is at most diffusive \cite{hsu2002stochastic} (i.e., has typical displacement $O(\sqrt{t})$ at time $t$), it is reasonable to conjecture that the same is true for random walk on non-negatively curved graphs of bounded-degree and that Salez's theorem can be improved from $o(t)$ to $O(\sqrt{t})$.

\medskip

This paper gives a complete affirmative answer to these questions for vertex-transitive graphs (i.e., those graphs for which the automorphism group acts transitively on the vertex set), with our results applying in particular to Cayley graphs of (finitely generated) groups.

\begin{theorem}\label{thm:growth_main}
	Let $G=(V,E)$ be a  transitive graph  of degree $d\in\mathbb{N}$. If $G$ has non-negative Ollivier--Ricci curvature
	then it has polynomial growth with dimension $O(d^3)$. More precisely, there exists a universal constant $C_0$ such that
	\begin{equation}\notag 
		\# B(x,r)\leq  r^{C_0 d^3}
	\end{equation}
  for every $x\in V$ and every $r\geq 2$.
\end{theorem}

\begin{corollary}
\label{cor:diffusivity}
For each integer $d\geq 1$ there exists a constant $C(d)<\infty$ such that if $G$ is a transitive graph of non-negative Ollivier--Ricci curvature and vertex degree $d$ then $\mathbf{E}_x d(X_0,X_n)^2 \leq C(d)n$ for every $n\geq 1$.
\end{corollary}

We have no reason to believe that the cubic dimension bound proven here is optimal. \cref{cor:diffusivity} follows from \cref{thm:growth_main} together with the finitary structure theory of transitive graphs of polynomial growth \cite{TesseraTointon,bgt12} 
and the proof does not yield any explicit control of the constant $C(d)$.

\medskip

While we are not yet able to prove polynomial growth for non-transitive bounded-degree graphs of non-negative Ollivier--Ricci curvature, we are able to prove a \emph{quasi-polynomial} $r^{O(\sqrt{\log r})}$ growth bound in a certain averaged sense, greatly strengthening the results of \cite{hutchcroft2024relation,salez2022sparse}:

\begin{thm}
\label{thm:main_finite}
There exists a universal constant $A$ such that the following holds.
 If $G=(V,E)$ is a finite, connected graph with non-negative Ollivier--Ricci curvature and degrees bounded by $d$ and degree ratios bounded by $C_\mathrm{rev}\defeq \max_{x,y\in V}\deg(x)/\deg(y)$, then
\begin{align*}
  \frac{1}{|V|}\sum_{x\in V} \log \# B(x,r) &\leq A d^3 C_\mathrm{rev}^{7/2}(\log r)^{3/2}
  \qquad \text{ and }\\
   \frac{1}{|V|}\sum_{x\in V} \mathbf{E}_x \!\left[d(X_0,X_n)^2\right] &\leq A d^3 C_\mathrm{rev}^{13/2} n (\log n)^{3/2}
\end{align*}
 for every $r,n\geq 2$. 
\end{thm}

Taking $r$ to be the diameter of the graph, we obtain that every finite, bounded-degree graph of non-negative Ollivier--Ricci curvature has volume at most $\exp\big(Ad^3 C_\mathrm{rev}^{7/2} (\log( \operatorname{diam}+1))^{3/2}\big)$ for some universal constant $A$.

\begin{remark}
It was recently proven by \cite{zhou2026pointwise} that every graph of non-negative Ollivier--Ricci curvature admits the significantly weaker but \emph{pointwise} bounds 
\begin{align*}
\# B(x,r) &\leq \exp\left[ \exp\left[ C(d) \sqrt{\log r \log\log r}\right]\right]
  \qquad \text{ and }\\
 \mathbf{E}_x \!\left[d(X_0,X_n)^2\right] &\leq n  \exp\left[ C(d) \sqrt{\log n \log\log n}\right]
\end{align*}
for all $r,n\geq 3$.
This result appeared before the present version of this paper but after the first version of this paper (in which a similar bound was proven in the same averaged sense considered in \cref{thm:main_finite}). We conjecture that bounds of the form
\begin{align*}
\# B(x,r) \leq C_d r^{C_d}
  \qquad \text{ and } \qquad 
 \mathbf{E}_x \!\left[d(X_0,X_n)^2\right] \leq C_d n
\end{align*}
hold pointwise for all graphs of non-negative Ollivier--Ricci curvature with degrees bounded by $d$.
\end{remark}

\begin{remark}
The bounds of \cref{thm:main_finite} also hold for \emph{unimodular random rooted graphs} of non-negative Ollivier--Ricci curvature and degrees bounded by $d$, with the same proof. Similarly, the proof of \cref{thm:growth_main} extends  to \emph{quasi-transitive} graphs (i.e., graphs for which the automorphism group has finitely many vertex orbits). 
\end{remark}

\noindent \textbf{Structural consequences.} It is a celebrated theorem of Gromov \cite{MR623534} (see also \cite{ozawa2018functional,kleiner2010new}) that a group admits a Cayley graph of polynomial volume growth if and only if it is \emph{virtually nilpotent}, meaning that it has a nilpotent subgroup of finite index.  \emph{Trofimov's theorem} \cite{MR735714} extends this structure theory to transitive graphs,  stating that every transitive graph of polynomial volume growth is quasi-isometric
 to a Cayley graph (which must be the Cayley graph of a virtually nilpotent group by Gromov's theorem). As such, it follows from \cref{thm:growth_main} that every finitely generated group admitting a Cayley graph of non-negative Ollivier--Ricci curvature must be virtually nilpotent and more generally that every locally finite transitive graph of non-negative Ollivier--Ricci curvature must be quasi-isometric to the Cayley graph of a nilpotent group.

In \cite{CamilloElia}, the first named author and E.\ Bruè proved that finitely generated, virtually nilpotent groups admitting Cayley graphs of non-negative Ollivier--Ricci curvature must, in fact, be virtually \emph{abelian}, and more generally that transitive graphs of polynomial growth and non-negative Ollivier--Ricci curvature must be quasi-isometric to $\Z^k$ for some $k\geq 0$. It is also shown conversely that every finitely generated virtually abelian group admits a Cayley graph of non-negative Ollivier--Ricci curvature. As such, \cref{thm:growth_main} has the following very strong structural consequences.

\begin{corollary}\label{cor:virtually_abelian}
A finitely generated group $G$ admits a Cayley graph of non-negative Ollivier--Ricci curvature if and only if it is virtually abelian.
  Quantitatively, there exists a universal constant $C_0$ such that if a Cayley graph $\operatorname{Cay}(G,S)$ has non-negative Ollivier--Ricci curvature then $G$ admits a finite-index subgroup isomorphic to the free abelian group $\Z^k$ for some $k\leq C_0 |S|^3$.
\end{corollary}

\begin{corollary}\label{cor:quasiZd}
Every locally finite transitive graph of non-negative Ollivier--Ricci curvature is quasi-isometric to a Euclidean lattice. Quantitatively,  there exists a universal constant $C_0$ such that if $G=(V,E)$ is a transitive graph of non-negative Ollivier--Ricci curvature with degrees bounded by $d$ then $G$ is quasi-isometric to $\Z^k$ for some $k\leq C_0 d^3$.
\end{corollary}

These corollaries can be thought of as a discrete analogue of the classical fact that fundamental groups of closed manifolds with non-negative Ricci curvature have virtually abelian fundamental groups (itself a simple consequence of the Cheeger--Gromoll splitting theorem \cite{cheeger1971splitting}).  Analogues of this theorem have previously been proven for groups satisfying stronger and more algebraic notions of non-negative curvature in \cite{nguyen2023cheeger,bar2022conjugation}.

While \cref{cor:virtually_abelian,cor:quasiZd} are vacuous for \emph{finite} graphs, it is possible to deduce the following strong finitary statement from \cref{cor:virtually_abelian} together with the fact that for each $d\geq 1$ there is a collection of finitely many degree-$d$ Cayley graphs of non-negative Ollivier--Ricci curvature having all other such graphs as quotients 
(\Cref{finitaryproof}). 

\begin{corollary}
\label{finitary}
For every $d\in\N$, there exists a constant $C_d\in\N$ such that the following holds.
    Let $G$ be a  finitely generated group and let $S$ be a finite, symmetric, set of generators with $|S|\le d$. If the Cayley graph $\operatorname{Cay}(G,S)$ has non-negative Ollivier--Ricci curvature, then $G$ contains an abelian subgroup of index bounded by $C_{d}$.
\end{corollary}

 We expect that one can also formulate and prove a similar finitary
version of \cref{cor:quasiZd}, but this appears to require one to first prove a finitary version of the results of \cite{CamilloElia} on \emph{large-scale Ollivier--Ricci curvature} and we do not pursue this here.

\begin{remark}[The Lin--Lu--Yau curvature and the role of laziness]
Let us briefly discuss the role of the laziness parameter. Generally, the normalized Ollivier curvature is increasing with the laziness parameter $\alpha$. However, for $\alpha \geq \frac 1 2$, the normalized curvature stays constant in $\alpha$, see \cite{loisel2014ricci}.
As such, the condition of non-negative Ollivier curvature is equivalent for all choices of the laziness parameter between one half and one.
As a consequence, our results also apply to graphs of non-negative \emph{Lin--Lu--Yau curvature} \cite{lin2011ricci} since this curvature is non-negative if and only if the Ollivier--Ricci curvature is non-negative \cite[Proposition 2]{loisel2014ricci}. 
\end{remark}

\begin{remark}[Other notions of curvature]
The property of having non-negative Ollivier--Ricci curvature is weaker than that of having non-negative \emph{sectional} curvature \cite[Section 4.2]{salez2025modern}, in which the $L^1$ Wasserstein metric is replaced by the $L^\infty$ Wasserstein metric or, equivalently, the condition that the random walks can be coupled so that $d(X_n,Y_n)$ is a supermartingale is replaced by the condition that they can be coupled so that $d(X_n,Y_n)$ is non-increasing. As such, all our results also apply automatically to graphs of non-negative sectional curvature. It was recently proven by Guo, Huang, and Huang \cite{guo2026nonnegative} (see also \cite{blachar2026edge}) that non-negative \emph{Bakry--Emery} curvature implies volume doubling and hence polynomial growth for all bounded-degree graphs; this is a different notion of curvature than Ollivier--Ricci and our result is neither weaker nor stronger than the results of \cite{guo2026nonnegative}.
\end{remark}

\medskip

\noindent \textbf{About the proof.} 
The starting point of our analysis is a well-known bound (arising from the coupling method) on the total variation distance between the step-$n$ random walk distributions starting from two neighbouring vertices in a graph of non-negative Ollivier--Ricci curvature:
\begin{equation}
  \operatorname{TV}_n \defeq \sup_{x\sim y} \|P^n(x,\cdot)-P^n(y,\cdot)\|_{\TV} \preceq \sqrt{\frac{d}{n}},\notag
\end{equation}
see e.g.\ \cite[Proposition 16]{salez2022sparse} for this inequality, \cite{bubley1997path} for the coupling method, and \cite[Corollary~3.1.2]{munch2023non} for a continuous time version. 
This bound is stated in detail in \cref{thm:Salez}.
As in \cite{hutchcroft2024relation}, we use the kernels $(P^n(x,\cdot))_{x\in V}$ to define a random partition of our vertex set into cells $\{[x]_n:x\in V\}$ with the average boundary-to-volume ratio of these cells governed by the quantity $\operatorname{TV}_n$ and with the volume and diameter of cells governed by the entropy and typical displacement of the random walk. A detailed description of this construction is given in \cref{sec:the_grand_poisson_cell_decomposition}.
The main technical innovations of our work compared to the earlier works \cite{salez2022sparse,hutchcroft2024relation} are as follows:
\begin{enumerate}
  \item We make crucial use of an isoperimetric \emph{lower bound} for non-negatively curved graphs due to the third named author \cite{munch2023ollivier}, which we generalize (\cref{thm:munch}) to establish a bound of the form
  \begin{equation}
  \label{eq:Munch_intro}
    \frac{|\partial W|}{|W|} \gtrsim  \frac{1}{\diam(\Lambda)} \log \frac{|\Lambda|}{|W|}
  \end{equation}
  for every two finite sets $W \subseteq \Lambda$ in a graph of non-negative Ollivier--Ricci curvature. (Here we have suppressed degree dependencies, see \cref{thm:munch} for the precise inequality.)
    To keep the paper self-contained, a full proof of this inequality  is given in \Cref{sec:Munch_appendix}.
  \item We prove in \cref{subsec:entropy_and_cell_size} that the cell sizes in the decomposition are \emph{equivalent} to the random walk entropy in the sense that $\E \log \# [x]_n \asymp H_n$ where $H_n$ is the Shannon entropy of the random walk at time $n$. (For finite graphs both quantities must be averaged over the initial vertex for this to hold.) This allows for a sharper and more efficient conversion of bounds on cell sizes into control of the random walk compared to \cite{hutchcroft2024relation}.  (In \cite{hutchcroft2024relation} it was observed only that cell sizes admit \emph{upper bounds} in terms of the entropy.) This two-sided relationship between cell sizes and entropy is crucial for the derivation of our functional inequalities.
  \item In \cref{subsec:Azuma} we observe that if $f:V\to \R$ is a $1$-Lipschitz function on a graph of non-negative Ollivier--Ricci curvature and $X$ is the random walk on $G$ started at some vertex $x$ then 
  \[(\mathbf{E}_x[f(X_n) \mid X_m])_{m=0}^n=(P^{n-m}f(X_m))_{m=0}^n\] defines a martingale with increments bounded by $2$, so that $f(X_n)-\mathbf{E}_x f(X_n)$ satisfies sub-Gaussian concentration estimates by Azuma--Hoeffding. These estimates (applied with $f=d(x,\cdot)$) are extremely useful when analyzing the geometry of cells in the random cell decomposition associated to the random walk. One particularly important consequence of this sub-Gaussian estimate is \cref{prop:displacement_doubling}, which bounds the displacement $D_n(x)=\mathbf{E}_x d(X_0,X_n)$ in terms of the volume doubling around $x$ on scale $D_n(x)$.
\end{enumerate}
In both the transitive and finite cases, we combine all of these ingredients into functional inequalities for the volume growth and expected log cell size which capture the tension between the isoperimetric lower bound \eqref{eq:Munch_intro} and the isoperimetric upper bounds attained through the random cell decomposition; all three ingredients we have just discussed contribute importantly to the form of these functional inequalities. Our bounds on the volume growth are then proven by an inductive analysis of these functional inequalities.

\medskip
\noindent \textbf{Notation.} We write $\asymp$, $\preceq$, and $\succeq$ for equalities and inequalities holding to within universal positive multiplicative constants (which may change at every occurrence). Given a set of vertices $W$, we write $\partial W$ for the set of edges with one endpoint in $W$ and the other in $V\setminus W$. To account for the possibility that our graphs have multiple edges, all sums of the form $\sum_{x\sim y}$ will be taken with multiplicity according to the number of edges connecting $x$ and $y$.

\section{Background and preliminary estimates}
\label{sec:background_and_preliminary_estimates}

\subsection{Total variation and the isoperimetric lower bound}

In this section we precisely state the total variation estimate and the isoperimetric lower bound generalizing the earlier work of the third named author.

\medskip

\noindent \textbf{The total variation estimate.}
We begin by stating the total variation estimate \cite[Proposition 16]{salez2022sparse} (as stated in \cite[Page 15]{munch2023mixing}), which will play a central role in our analysis.
Recall that the \textbf{total variation distance} between the probability measures $\mu$ and $\nu$ on a countable set S
is defined to be
\begin{align*}
  \|\mu-\nu\|_\TV &\defeq \frac{1}{2}\sum_{x\in S} |\mu(x)-\nu(x)| \\&= \inf \left\{\P(X \neq Y) : \P \text{ a probability measure on $S^2$ with marginals $\mu$ and $\nu$}\right\};
\end{align*}
see e.g.\ \cite[Proposition 4.7]{MCMT} for a proof that these two quantities are equivalent.

\begin{theorem}
\label{thm:Salez}
Let $G=(V,E)$ be a   graph with non-negative Ollivier--Ricci curvature and with degrees bounded by $d$. Then
\begin{equation}\notag
  \|P^n(x,\cdot)-P^n(y,\cdot)\|_{\TV} \leq \sqrt{\frac{20 d}{n+1}}\qquad\text{for every $x\sim y$ and $n\in\mathbb{N}$}.
\end{equation}
\end{theorem}

An important consequence of this theorem is that bounded-degree graphs of non-negative Ollivier--Ricci curvature are \emph{amenable}, meaning that there exist finite sets of vertices with arbitrarily small boundary-to-volume ratio $|\partial W|/|W|$; see \cite{hutchcroft2024relation} for details.

\begin{corollary}
\label{cor:amenable}
Every bounded-degree graph with non-negative Ollivier--Ricci curvature is amenable.
\end{corollary}

 (For groups, the fact that \cref{thm:Salez} implies amenability is an immediate consequence of Reiter's condition for amenability with $(P^n(\mathrm{id},\cdot))_{n\geq 1}$ being a sequence of Reiter functions. Note also that \cref{thm:Salez} easily implies that bounded-degree graphs of non-negative Ollivier--Ricci curvature are Liouville, which implies amenability in the transitive case. This is the only setting \cref{cor:amenable} will be applied in.)

\medskip

\noindent\textbf{The isoperimetric lower bound.}
As mentioned in the introduction, the proof of our main theorems will exploit the tension between the isoperimetric \emph{upper bounds} derived using the grand Poisson cell decomposition (introduced in \cref{sec:the_grand_poisson_cell_decomposition}) with the isoperimetric \emph{lower bounds} of \cite{munch2023ollivier}, which we now state. We write $\pi(W)=\sum_{x\in W} \deg(x)$ and $\diam(W)=\max_{x,y\in W}d(x,y)$ for the stationary measure and diameter of the set $W \subseteq V$ and write $\partial W$ for the set of \emph{edges} with one endpoint in $W$ and the other in $V\setminus W$.

\begin{theorem}
\label{thm:munch}
If $G=(V,E)$ is a connected graph with non-negative Ollivier--Ricci curvature and maximum degree at most $d$ then
\[
  \frac{\#\partial W}{\pi(W)} \geq \frac{1}{4d\diam(\Lambda)}\log \left[\frac{\pi(\Lambda)}{2\pi(W)}\right] 
\]
for every pair of finite sets $\emptyset\ne W\subseteq \Lambda \subseteq V$.
\end{theorem}

This bound is an easy generalization of \cite[Theorem 6.1]{munch2023ollivier} (in which $G$ is finite and $\Lambda$ is taken to be the entire vertex set). A full proof is given in \Cref{sec:Munch_appendix}. 
For transitive graphs (or, more generally, regular graphs) this bound can be written as stating that
\begin{equation}\label{sdcdsacsc}
  \log \#\Lambda\leq \log \#W+ 4\diam(\Lambda)\frac{\#\partial W}{\#W}+\log 2
\end{equation}
whenever $W\subseteq \Lambda$ are finite sets of vertices.

\subsection{Unimodularity and the mass-transport principle}

Let $G=(V,E)$ be a locally finite  transitive graph. We say that $G$ is \textbf{unimodular} if it satisfies the \textbf{mass-transport principle}, meaning that if 
 $F:V\times V\rightarrow [0,\infty]$ which is invariant under the diagonal action of $\Aut(G)$ in the sense that
  $F(\gamma x,\gamma y)=F(x,y)$ for every $\gamma\in\Aut(G)$ and $x,y\in V$
then
\begin{equation}\label{MTP}
  \sum_{x\in V}F(o,x)=\sum_{x\in V}F(x,o).
\end{equation}
for every vertex $o\in V$. Every Cayley graph is unimodular, with the mass-transport principle being a trivial consequence of the fact that left multiplication defines a graph automorphism and inversion is bijective:
  \begin{equation}\notag
    \sum_{g\in V}F(\mathrm{id},g)=\sum_{g\in V}F(g^{-1}, \mathrm{id})= \sum_{g\in V}F(g, \mathrm{id}).
  \end{equation}
 The mass-transport principle is also trivially verified for every \emph{finite} transitive graph. 
It is a classical theorem of Soardi and Woess \cite{MR1082868} that every amenable transitive graph is unimodular, and since every non-negatively curved transitive graph is amenable by \cref{cor:amenable} it follows that every non-negatively curved transitive graph is unimodular also.

\subsection{Concentration of Lipschitz functions of the random walk}
\label{subsec:Azuma}

In this section we prove the following proposition establishing sub-Gaussian estimates on $f(X_n)$ when $X$ is a random walk on a graph of non-negative Ollivier--Ricci curvature and $f$ is a Lipschitz function on $G$. Note that the constants in this estimate are independent of the degree. In addition to various other applications, this estimate will be used to bound the expected displacement in terms of the local volume doubling dimension in the next section.

\begin{prop}
  \label{prop:Lipschitz_concentration}
  Let $G=(V,E)$ be a locally finite graph of non-negative Ollivier--Ricci curvature.
  If $f:V\to \R$ is Lipschitz then, for every $\lambda>0$,
  \begin{equation}\notag
  \mathbf{P}_x\left( |f(X_n)-\bE_x f(X_n)| \geq \lambda n^{1/2} \Lip(f) \right) \leq 2\exp\Big(-\frac{1}{8}\lambda^2\Big)\quad\text{  for every $x\in V$ and $n\geq 1$}.
  \end{equation}
\end{prop}

\begin{proof}[Proof of \Cref{prop:Lipschitz_concentration}]
  It suffices by scaling to consider the case $\Lip(f)=1$ (the case $\Lip(f)=0$ being trivial). Fix $n\geq 1$ and consider the martingale 
  \begin{align*}
    M_m &= \mathbf{E}_x[f(X_n) \mid (X_i)_{i=0}^m ]-\mathbf{E}_x f(X_n) = [P^{n-m}f](X_m)-[P^nf](x),
  \end{align*}
  which has $M_0=0$ and $M_n = f(X_n)-\mathbf{E}_x f(X_n)$. The increments of this martingale satisfy
  \begin{multline*}
    M_{m+1}-M_m= [P^{n-m-1} f](X_{m+1})-[P^{n-m} f](X_{m}) \\=[P^{n-m-1} f](X_{m+1})-[P^{n-m} f](X_{m+1}) + [P^{n-m} f](X_{m+1}) -[P^{n-m} f](X_{m}).
  \end{multline*}
  Since $G$ has non-negative Ollivier--Ricci curvature, $P$ is a contraction on the space of Lipschitz functions. This implies that $|[P^{n-m} f](X_{m+1}) -[P^{n-m} f](X_{m})|\leq 1$, while we also have trivially that
  \[
  |[P^{n-m-1} f](w)-[P^{n-m} f](w)| = |\mathbf{E}_w f(X_{n-m-1})-\mathbf{E}_w f(X_{n-m})| \leq \mathbf{E}_w |f(X_{n-m-1})- f(X_{n-m})| \leq 1.
  \]
  It follows that the martingale $M$ has increments bounded by $2$, and the claim follows from the Azuma--Hoeffding inequality.
\end{proof}

\begin{remark}
  Antoine Song has pointed out to us that this proposition can be viewed as a generalization of the Maurey--Pisier theorem \cite[p.\ 181]{pisier2006probabilistic} (which establishes sub-Gaussian estimates on Lipschitz functions of Gaussians in $\R^d$) to the non-negatively curved context.
\end{remark}

\subsection{Diffusivity from doubling}

In this section we apply \cref{prop:Lipschitz_concentration} to prove the following proposition, which bounds
the \textbf{expected displacement} of the random walk $D_n(x)\defeq\E_x d(x,X_n)$ from above in terms of the \emph{local volume doubling dimension}.

\begin{prop}\label{prop:displacement_doubling}
Let $G$ be a locally finite graph of non-negative Ollivier--Ricci curvature. Then
\[
D_n(x) \leq \sqrt{ 32 n \log \frac{9\pi(B(x,4D_n(x)))}{2\pi(B(x,D_n(x)/8))}}
\]
for every $x\in V$ and $n\geq 1$.
\end{prop}

The proof of this proposition will use the following lemma.

\begin{lemma}
\label{lem:D_n_Lipschitz}
If $G=(V,E)$ has non-negative Ollivier--Ricci curvature then the function $x\mapsto D_n(x)$ is $2$-Lipschitz on $V$.
\end{lemma}

\begin{proof}[Proof of \cref{lem:D_n_Lipschitz}]
Fix two neighbouring vertices $x$ and $y$. Writing $d_x$ and $d_y$ for the functions $d_x(z)=d(x,z)$ and $d_y(z)=d(y,z)$, we can write
\[
  |D_n(x)-D_n(y)|=|(P^n d_x)(x)-(P^nd_y)(y)| \leq |(P^n d_x)(x)-(P^nd_x)(y)| + |(P^n d_x)(y)-(P^nd_y)(y)|.
\]
The first term on the right hand side is at most $1$ since $\Lip(P^n d_x)\leq \Lip(d_x) =1$ by the defining property of non-negative Ollivier--Ricci curvature, while the second is equal to $|\mathbf{E}_y [d(x,X_n)-d(y,X_n)]|$ which is also trivially at most one.
\end{proof}

The proof of \cref{prop:displacement_doubling} will also rely on the fact that lazy random walk on any locally finite graph satisfies
\begin{multline}\label{vedscdasc}
    \mathbf{P}_A(X_n \in B) =\frac{1}{\pi(A)}\sum_{a\in A} \deg(a) \mathbf{P}_a (X_n \in B) = \frac{1}{\pi(A)}\langle P^n \mathbbm{1}_B,\mathbbm{1}_A\rangle_\pi \\\leq \frac{1}{\pi(A)}
    \sqrt{\langle P^n \mathbbm{1}_A,\mathbbm{1}_A\rangle_\pi\langle P^n \mathbbm{1}_B,\mathbbm{1}_B\rangle_\pi}
    =\sqrt{\frac{\pi(B)}{\pi(A)} \mathbf{P}_{A}(X_n \in A)\mathbf{P}_{B}(X_n \in B)}
\end{multline}
for any two finite sets $A$ and $B$, where $\mathbf{P}_A$ denotes the law of a random walk started from a point of $A$ chosen with weight proportional to its degree. The inequality appearing here follows from the fact that the lazy random walk $n$-step transition matrix $P^n$ is self-adjoint and positive semidefinite with respect to the inner product $\langle f,g\rangle_\pi\defeq\sum_{x\in V} \deg(x)f(x)g(x)$ so that $\langle P^n \cdot, \cdot\rangle_\pi$ is an inner product and hence satisfies the Cauchy--Schwarz inequality. (The fact that $P^n$ is positive semidefinite follows by writing $P=\frac{1}{2}(I+Q)$ with $Q$ the non-lazy simple random walk transition matrix, which has spectrum in $[-1,1]$.)

\begin{proof}[Proof of \cref{prop:displacement_doubling}]
%
%
Fix $x\in V$ and $n\geq 1$ and let $B_1$ and $B_2$ be the balls of radius $D_n(x)/8$ and $4D_n(x)$ around $x$ respectively. 
Since $y\mapsto D_n(y)$ is $2$-Lipschitz by \cref{lem:D_n_Lipschitz}, every $y\in B_1$ has $\frac{3}{4}D_n(x) \leq D_n(y) \leq \frac{5}{4}D_n(x)$ and hence 
\begin{equation}\notag
  \mathbf{P}_y (X_n \in B_2) \geq 
  \mathbf{P}_y\left(d(y,X_n) \leq \frac{31}{8} D_n(x)\right) \geq \mathbf{P}_y\left(d(y,X_n) \leq \frac{31}{10} D_n(y)\right) 
   \geq \frac{2}{3}
\end{equation}
by Markov's inequality. Since this estimate holds for every $y\in B_1$, it follows in particular that
 \begin{equation}\label{vrecdscadscasc}
  \mathbf{P}_{B_1}(X_n \in B_2) \geq \frac{2}{3}.
 \end{equation}
On the other hand, \Cref{prop:Lipschitz_concentration} applied to $d(y,\cdot)$ yields that
\begin{multline*}
\mathbf{P}_{B_1}(X_n \in B_1) \leq \max_{y\in B_1} \mathbf{P}_y\left(d(y,X_n) \leq \frac{1}{4}D_n(x)\right) \\\leq \max_{y\in B_1} 
\mathbf{P}_y\left(|d(y,X_n)-D_n(y)| \geq \frac{1}{2}D_n(x)\right)  \leq 2 e^{-\frac{D_n(x)^2}{32 n}}.
\end{multline*}
Substituting this estimate and \eqref{vrecdscadscasc} into \eqref{vedscdasc} and using the trivial bound $\mathbf{P}_{{B_2}}(X_n \in B_2)\leq 1$ yields that
\[
\frac{2}{3} \leq \sqrt{\frac{\pi(B_2)}{\pi(B_1)}2 e^{-\frac{D_n(x)^2}{32n}}}
\]
which rearranges to yield the claim.
\end{proof}

\section{The grand Poisson cell decomposition}
\label{sec:the_grand_poisson_cell_decomposition}

In this section we recall the definition of the random cell decomposition used in \cite{hutchcroft2024relation}, which we refer to as the \emph{grand Poisson cell decomposition} associated to the random walk kernels $(P^n(x,\cdot))_{x\in V}$, and prove a relationship between the volume of cells in this decomposition and the Shannon entropy of the random walk. We begin by discussing (what we call) the \emph{grand Poisson coupling} \cite{angel2019pairwise}, a simultaneous coupling of all probability measures on a countable set that attains the total variation distance between any two such measures up to a factor of $2$.

\medskip

\noindent
\textbf{The grand Poisson coupling.}
Let $S$ be a countable set. We consider the following simultaneous coupling of all probability measures on $S$:
Let $\mathcal{P}$ be a Poisson point process on $S\times [0,1]\times [0,\infty)$ with intensity the product of the counting measure with two copies of the Lebesgue measure, and let $\ell:S\times [0,1]\times [0,\infty)\to [0,\infty)$ be the projection to the last coordinate, which we call the \textbf{label} of the point. For each probability measure $\mu$ on $S$, let 
\[
  X(\mu) = \argmin_x \min \Bigl\{\ell(\xi) : \xi \in \mathcal{P} \cap (\{x\}\times [0,\mu(x)] \times [0,\infty)) \Bigr\}.
\]
(Note that $X(\mu)$ is well defined for every $\mu$ on the probability-one event that no two points in the Poisson point process have the same label.)
We call this the \textbf{grand Poisson coupling} on $S$. Equivalently, one can think of points appearing on $S\times [0,1]$ continuously in time, where the rate at which points appear in a set $A \subseteq S\times [0,1]$ is proportional to the measure of $A$, the arrival times of the points being the labels in the previous description. In this picture, the random variable $X(\mu)$ is defined to be the $S$-coordinate of the first point that appears in the set $\{(x,p):x\in S, 0\leq p \leq \mu(x)\}$.

One can easily verify that $\P(X(\mu)=x)=\mu(x)$,
so that the grand Poisson coupling does indeed define a simultaneous coupling of all probability measures on $S$.
We write $\Xi(\mu)$ for the point of $\mathcal{P}$ determining $X(\mu)$, so that $X(\mu)$ is the $S$-coordinate of $\Xi(\mu)$ and in particular $X(\mu)=X(\nu)$ whenever $\Xi(\mu)=\Xi(\nu)$.
 It is observed  in \cite{angel2019pairwise}  that
\begin{equation}
\label{eq:coupling_TV}
  \P(X(\mu)\neq X(\nu)) \leq \P(\Xi(\mu)\neq \Xi(\nu)) = \frac{1}{1+\|\mu-\nu\|_\TV} \sum_{x\in S} |\mu(x)-\nu(x)| = \frac{2 \|\mu-\nu\|_\TV}{1+\|\mu-\nu\|_\TV}.
\end{equation}
 Thus, the grand Poisson coupling defines a simultaneous coupling of all probability measures on $S$ that is optimal up to a factor of $2$ for each pair of measures $\mu$ and $\nu$. This identity has the simple interpretation that $\Xi(\mu)=\Xi(\nu)$ precisely when the first point that appears in the set $\{(x,p): x\in S, p \leq \mu(x)\vee \nu(x)\}$ belongs to both sets $\{(x,p): x\in S, p \leq \mu(x)\}$ and $\{(x,p): x\in S, p \leq \nu(x)\}$. This is equivalent to a point appearing in the set $\{(x,p):x\in S, p\leq \mu(x)\wedge \nu(x)\}$ before a point appears in $\{(x,p): x\in S, \mu(x) \wedge \nu(x) \leq p\leq \mu(x)\vee \nu(x)\}$. This happens with probability
\[
 \P(\Xi(\mu)= \Xi(\nu)) = \frac{\sum_{x\in S} \mu(x) \wedge \nu(x)}{\sum_{x\in S} \mu(x) \vee \nu(x)}=\frac{1-\|\mu-\nu\|_\TV}{1+\|\mu-\nu\|_\TV},
\]
where we used that
$x \vee y + x \wedge y = x+y$ and 
 $x \vee y - x \wedge y = |x-y|$
and hence that
\[
x\vee y = \frac{x+y+|x-y|}{2} \qquad \text{ and } \qquad
x\wedge y = \frac{x+y-|x-y|}{2}
\]
in the last equality.

\medskip

We will need the following lemma, which follows by the same reasoning used to derive \eqref{eq:coupling_TV}.

\begin{lemma}
\label{lem:coupling_two_point}
Consider the grand Poisson coupling on a countable set $S$.
If $\mu$ and $\nu$ are two probability measures on the countable set $S$ then
\[
  \frac{\mu(x) \wedge \nu(x)}{1+\|\mu-\nu\|_\TV} \leq \P(X(\mu)=X(\nu)=x) \leq \mu(x) \wedge \nu(x)
\]
for every $x\in S$.
Equivalently,
\[
  \frac{1}{1+\|\mu-\nu\|_\TV} \left[\frac{\nu(x)}{\mu(x)}\wedge 1\right] \leq \P(X(\nu)=x \mid X(\mu)=x) \leq \frac{\nu(x)}{\mu(x)}\wedge 1
\]
for every $x\in S$ with $\mu(x)\neq 0$.
\end{lemma}


\begin{proof}[Proof of \Cref{lem:coupling_two_point}]
The upper bound is immediate from the fact that $\P(X(\mu)=X(\nu)=x)$ is at most $\P(X(\mu)=x)\wedge \P(X(\nu)=x)=\mu(x)\wedge \nu(x)$; we focus on the lower bound.
Fix $x\in S$, write $\TV=\|\mu-\nu\|_\TV$ and assume without loss of generality that $\mu(x) \geq \nu(x)$.
Let $\tau$ be the label of the point determining the value of $X(\mu)$. The value of $\tau$ is independent of $X(\mu)$.
Given $X(\mu)=x$ and $\tau=t$, the probability that the point with label $\tau$ also sits under the graph of $\nu$ is $\nu(x)/\mu(x)$.
Given that this occurs, $\Xi(\nu)=\Xi(\mu)$ if and only if there are no points of the Poisson point process in $\{(w,p,s): w\in S, \mu(w) \leq p \leq \nu(w), s \leq t\}$, this set being the part of the region lying under the graph of $\nu$ but not under the graph of $\mu$. Since the restrictions of the Poisson point process to disjoint sets are independent and the volume of the set $\{(w,p): w\in S, \mu(w) \leq p \leq \nu(w)\}$ is equal to $\TV$, the conditional probability that this occurs 
 is $e^{-t\operatorname{TV}}$. Since $X(\nu)=X(\mu)=x$ whenever this occurs, the probability that $X(\mu)=X(\nu)=x$ is at least
\[
\mu(x)\cdot \frac{\nu(x)}{\mu(x)} \int_0^\infty e^{-t \operatorname{TV}} e^{-t} \dif t = \frac{\nu(x)}{1+\operatorname{TV}} 
\]
as claimed. 
\end{proof}

\noindent\textbf{The grand Poisson cell decomposition.}
Now suppose that we have a kernel $\kappa:S^2\to [0,1]$ on a countable set $S$ with the property that $\sum_{y\in S}\kappa(x,y)=1$ for every $x\in S$, so that $\mu_x(\{\,\cdot\,\})\defeq\kappa(x,\,\cdot\,)$ defines a probability measure on $S$ for each $x\in S$. We call such a function $\kappa:S^2\to[0,1]$ a \textbf{probability kernel} on $S$. (We will be mostly interested in the case that $\kappa(x,y)=P^n(x,y)$ is the $n$-step transition kernel of the lazy random walk on a graph, but formulate things more generally for now.) Considering the coupling of these measures induced by the grand Poisson coupling, we obtain a family of random variables $(X(x))_{x\in S}$ indexed by $S$. We use this family to define a partition of $S$ by setting the \textbf{cell} of $x$ to be $[x]=\{y\in S:X(y)=X(x)\}$. We call this the \textbf{grand Poisson cell decomposition} associated to the probability kernel $\kappa$.

\medskip

 Given a probability kernel $\kappa:V^2\to[0,1]$ we write $\mu_x$ for the probability measure with probability mass function $\mu_x(y)=\kappa(x,y)$ and say that $\kappa$ is \textbf{automorphism-invariant} if $\kappa(\gamma x, \gamma y)=\kappa(x,y)$ for every $x,y\in V$ and $\gamma \in \Aut(G)$.
Observe that if $\kappa:V^2\to[0,1]$ is automorphism-invariant then the resulting grand Poisson cell decomposition is also automorphism-invariant in the sense that if $\Pi$ is the partition of $V$ into the cells of the grand Poisson cell decomposition and $\gamma \in \Aut(G)$ then 
$\gamma \Pi \defeq \{\gamma C : C \in \Pi\}$
has the same law as $\Pi$. Indeed, if $(X(x))_{x\in V}$ are the random variables used to define the grand Poisson cell decomposition then the random variables $(Y(x))_{x\in V}$ defined by 
\begin{equation}\notag
Y(x)\defeq \gamma X(\gamma^{-1}x)\qquad\text{i.e.,} \qquad  Y(\gamma x)=\gamma X(x),
\end{equation}
have the same joint law as $(X(x))_{x\in V}$ by our  automorphism-invariance assumption on the probability kernel $\kappa$ defining $(\mu_x)_{x\in V}$ and $x$ and $y$ are in the same cell of $\gamma \Pi$ if and only if $Y(x)=Y(y)$. Indeed, 
if $\mathcal{P}$ is the Poisson point process on $V\times [0,1]\times [0,\infty)$ used to define the random variables $(X(x))_{x\in V}$, then  the random variables $(Y(x))_{x\in V}$ correspond to the grand Poisson decomposition with respect to the Poisson point process obtained by $\mathcal{P}$ through the map $(x,s,t)\mapsto (\gamma x,s,t)$.
In general we refer to any random partition $\Pi$ of our vertex set for which $\gamma\Pi$ has the same distribution as $\Pi$ for every $\gamma\in \Aut(G)$ as an \textbf{automorphism-invariant random partition} of $G$.

\medskip

On a unimodular transitive graph, the expected boundary-to-volume ratio of the cell is governed by the average total variation distance between neighbours as described in the following lemma. We say that a probability kernel $\kappa:V^2\to [0,1]$ on a graph has \textbf{finite range} if there exists $n<\infty$ such that $\kappa(x,y)=0$ for all $x,y\in V$ with $d(x,y) > n$; this forces all cells in the grand Poisson cell decomposition associated to $\kappa$ to have diameter at most $2n$ and hence to be finite.

\begin{lemma}
\label{lem:cell_isoperimetry_transitive}
Let $G=(V,E)$ be a locally finite unimodular transitive graph. If
$\Pi$ is an automorphism-invariant random partition of $V$ with finite cells then
\[
  \E \frac{\#\partial [o]}{\#[o]} = \sum_{x\sim o} \P(x\notin [o]),
\]
where we write $[o]$ for the cell of $o$ in $\Pi$. In particular, if
 $\kappa:V\times V\to [0,1]$ is an automorphism-invariant, finite-range probability kernel on $G$ then the grand Poisson cell decomposition associated to $\kappa$ satisfies
\begin{equation}\notag
	 \E \frac{\#\partial [o]}{\#[o]} =\sum_{z\sim o}\P(X(z)\ne X(o)) \leq \sum_{x\sim o} \frac{2\|\mu_x-\mu_o\|_\TV}{1+\|\mu_x-\mu_o\|_\TV}.
\end{equation}
\end{lemma}

\begin{proof}[Proof of \cref{lem:cell_isoperimetry_transitive}]
We apply the mass-transport principle $\sum_{x\in V}F(o,x)=\sum_{x\in V}F(x,o)$ to $F$ the function
\begin{equation}\notag
F(x,y)\defeq  \E \sum_{z\sim x}\mathbbm{1}([z]\ne [x])\frac{\mathbbm{1}([x]=[y])}{\#[x]},
\end{equation}
which is easily seen to be diagonally-invariant when $\Pi$ is an automorphism-invariant random partition.
The left-hand side is equal to
\begin{align*}
		\sum_{x\in V}F(o,x)&=\E \sum_{x\in V}\sum_{z\sim o} \mathbbm{1}([z]\ne  [o])\frac{\mathbbm{1}([o]=[x])}{\#[o]}
    =\E \sum_{z\sim o}\mathbbm{1}([z] \ne [o])
\end{align*}
while the right-hand side is equal to
\begin{equation*}
	\sum_{x\in V}F(x,o)=\E 	\sum_{x\in V}\sum_{z\sim x}\mathbbm{1}([z]\ne [x])\frac{\mathbbm{1}([x]=[o])}{\#[x]}= \E \sum_{x\in [o]}\sum_{z\sim x}\frac{\mathbbm{1}(z\notin [o])}{\#[o]}= \E \frac{\#\partial [o]}{\#[o]} .
\end{equation*}
The second claim follows from the first claim and \eqref{eq:coupling_TV}.
\end{proof}

We will also make use of the following lemma describing the effect of changing the probability kernel on the grand Poisson cell decomposition.

\begin{lemma}
	\label{lem:different_cells}
	Let $G=(V,E)$ be a locally finite unimodular transitive graph and
  suppose that $\kappa^1,\kappa^2:V\times V\to [0,1]$ are two automorphism-invariant finite-range probability kernels on $G$.
  If we fix $o\in V$ and let $[o]_1$ and $[o]_2$ be the cells of $o$ in the grand Poisson cell decompositions defined from $\kappa^1$ and $\kappa^2$ using the same Poisson point process $\cP$ then
	\begin{equation}\notag
		\E\frac{\#\big([o]_1\setminus [o]_2\big)}{\#[o]_1}\leq 4\|\mu_o^1-\mu_o^2\|_\TV.
	\end{equation}
\end{lemma}

\begin{proof}[Proof of \cref{lem:different_cells}]
Let $\mathcal{P}$ be the underlying Poisson point process on $V\times [0,1]\times [0,\infty)$ used to define both cell decompositions and for each $i=1,2$ let  $X^i=(X^i(x))_{x\in V}$ denote the associated family of $V$-valued random variables with laws $(\mu^i_x)_{x\in V}$.  As in the proof of \Cref{lem:cell_isoperimetry_transitive}, the function
\begin{equation}\notag
	 F(x,y)\defeq \E \left[\mathbbm{1}(X^1(y)\ne X^2(y))\frac{\mathbbm{1}(X^1(y)=X^1(x))}{\#[x]_1}\,\right],
\end{equation}
is invariant under the diagonal action of $\Aut(G)$
so that we can apply the mass-transport principle \eqref{MTP} to obtain that
	$\sum_{x\in V}F(o,x)=\sum_{x\in V}F(x,o)$.
The left-hand side of this equality is equal to
\begin{align*}
\E \sum_{x\in V}\mathbbm{1}(X^1(x)\ne X^2(x))\frac{\mathbbm{1}(X^1(x)=X^1(o))}{\#[o]_1}=\E\frac{1}{\#[o]_1}\sum_{x\in [o]_1}\mathbbm{1}(X^1(x)\ne X^2(x)),
\end{align*}
whereas the right-hand side is equal to
\begin{multline}\notag
		\E \sum_{x\in V} \mathbbm{1}(X^1(o)\ne X^2(o))\frac{\mathbbm{1}(X^1(o)=X^1(x))}{\#[x]_1}=\E  \mathbbm{1}(X^1(o)\ne X^2(o)) \sum_{x\in V}\frac{\mathbbm{1}(x\in [o]_1)}{\#[o]_1}\\=\P (X^1(o)\ne X^2(o))
\end{multline}
so that
\begin{equation}\label{ppefvd}
\E\frac{1}{\#[o]_1}\sum_{x\in [o]_1}\mathbbm{1}(X^1(x)\ne X^2(x))=\P (X^1(o)\ne X^2(o)).
\end{equation}
Note that if $X^1(o)=X^2(o)$ and $x\in [o]_1$ is such that $X^1(x)=X^2(x)$ then $x\in [o]_2$, so that
\begin{align*}
		\E\frac{\#\big([o]_1\setminus [o]_2\big)}{\#[o]_1} = \E \frac{1}{\#[o]_1}\sum_{x\in [o]_1}\mathbbm{1}(x\notin [o]_2)\leq \P(X^1(o)\neq X^2(o)) + \E \frac{1}{\#[o]_1}\sum_{x\in [o]_1}\mathbbm{1}(X^1(x)\ne X^2(x)).
\end{align*}
Using \eqref{ppefvd} we obtain that
\begin{align*}
		\E\frac{\#\big([o]_1\setminus [o]_2\big)}{\#[o]_1}&\leq \E \mathbbm{1}(X^1(o)\ne X^2(o))+\E   \frac{1}{\#[o]_1}\sum_{x\in [o]_1}\mathbbm{1}(X^1(x)\ne X^2(x))= 2\P (X^1(o)\ne X^2(o))\\
		&\leq  4 \|\mu_o^1-\mu_o^2\|_{\TV},
\end{align*}
where the last inequality follows from  \eqref{eq:coupling_TV}.
\end{proof}

We also have the following finite graph version of \cref{lem:cell_isoperimetry_transitive}
 that does not require any symmetry but instead requires us to average over the choice of the root vertex. For technical reasons we also extend the boundary identity to include a non-negative weight function.

\begin{lemma}
\label{lem:cell_isoperimetry_finite}
Let $G=(V,E)$ be a finite graph, let $\Pi$ be a random partition of $V$, and let $W:\Pi \to [0,\infty]$ be a (possibly random) function assigning a non-negative number to each cell of $\Pi$. Then we have that
\[
  \frac{1}{|V|}\sum_{x\in V} \E \left[W([x]) \frac{\#\partial [x]}{\#[x]}\right] = \frac{1}{|V|}\sum_{x\in V} \sum_{z\sim x} \E\left[ W([x]) \mathbbm{1}(z\notin [x])\right]
\]
where $[x]$ is the cell of $x$ in $\Pi$.
\end{lemma}


\begin{proof}[Proof of \cref{lem:cell_isoperimetry_finite}]
This proof is identical to that of \cref{lem:cell_isoperimetry_transitive} except that we use the trivial mass-transport principle for finite graphs 
\begin{equation}
  \frac{1}{|V|} \sum_{o\in V}\sum_{x\in V} F(o,x)= \frac{1}{|V|} \sum_{o\in V}\sum_{x\in V} F(x,o)
  \label{eq:finiteMTP}
\end{equation}
which holds for any function $F:V\times V\to [0,\infty]$, and use the weighted function
\[
  F(x,y)=\E \left[W([x])\sum_{z\sim x}\mathbbm{1}([z]\ne[x])\frac{\mathbbm{1}([x]=[y])}{\#[x]}\right]
\]
in the boundary computation.
\end{proof}

\subsection{Entropy and cell size}
\label{subsec:entropy_and_cell_size}

 We now relate the \emph{sizes} of the cells in the grand Poisson cell decomposition to the \emph{entropies} of the measures $(\mu_x)_{x\in S}$.
We will be particularly interested in the case that the kernels $\kappa$ arise from the (lazy) random walk on a bounded-degree graph. The relevant property of such kernels we will need is captured abstractly by the following definition: We define a kernel $\kappa$ to be \textbf{$C$-reversible} if
\[
C^{-1}\kappa(x,y)\leq \kappa(y,x)\leq C \kappa(x,y)
\]
for every $x,y\in S$. If $P^n$ is the $n$-step transition kernel of the (lazy) random walk on a bounded-degree graph,  $P^n$ is $C$-reversible with $C= \max_{x,y} \deg(x)/\deg(y)$.

\begin{prop}[Entropy is equivalent to cell size]
\label{prop:entropy_and_cell_size}
Consider the grand Poisson cell decomposition on a finite set $S$ associated to the probability kernel $\kappa:S^2\to [0,1]$. If $\kappa$ is $C$-reversible  for some $C\geq 1$ and we define $\mu_x(\,\cdot\,)=\kappa(x,\,\cdot\,)$ for each $x\in S$ then 
\[
 \frac{1}{8C^3|S|}\sum_{x\in S} \operatorname{Ent}(\mu_x) - \frac{1}{2}\leq \frac{1}{|S|}\sum_{x\in S}\E\log \#[x] \leq \frac{C}{|S|}\sum_{x\in S} \operatorname{Ent}(\mu_x)+2\log C.
\] 
\end{prop}

Here we recall that the Shannon entropy of a probability measure $\mu$ on a countable set $S$ is defined by $\operatorname{Ent}(\mu)=-\sum_x \mu(x)\log \mu(x)$.
We begin with the following basic lemma about the entropy of a single probability measure; we will apply this lemma with $\eps=1/2$.

\begin{lemma}
\label{lem:entropy_alternative}
For each $0<\eps\leq 1$ there exists a constant $C(\eps)$ such that
\[
(1-\eps) \operatorname{Ent}(\mu) - C(\eps) \leq 
\sum_x \mu(x) \log \#\{y:\mu(y)\geq \mu(x)\}\leq 
\sum_x \mu(x) \log \left[\sum_y \frac{\mu(y)}{\mu(x)} \wedge 1\right] \leq \operatorname{Ent}(\mu)
\]
for every probability measure $\mu$ on a countable set $S$. Moreover, for $\eps=1/2$ we can take $C(\eps)=1/4$.
\end{lemma}

\begin{proof}[Proof of \cref{lem:entropy_alternative}]
We have trivially that
\[\sum_x \mu(x) \log \left[\sum_y \frac{\mu(y)}{\mu(x)} \wedge 1\right]\leq \sum_x \mu(x) \log \left[\sum_y \frac{\mu(y)}{\mu(x)}\right] = \sum_x \mu(x) \log  \frac{1}{\mu(x)} =\operatorname{Ent}(\mu).\]
The middle inequality of the claim holds since each $y$ with $\mu(y)\geq \mu(x)$ contributes $1$ to the sum inside the logarithm. For the remaining lower bound, order the set $S=\{x_1,x_2,\ldots\}$ in decreasing order of $\mu(x_i)$ and let $\mu_i=\mu(x_i)$, so that $\#\{y : \mu(y)\geq \mu(x_k)\}\geq k$ for every $k\geq 1$ and hence
\[
 \sum_x \mu(x) \log \#\{y : \mu(y) \geq \mu(x)\} \geq \sum_{k=1}^\infty \mu_k \log k.
\]
Now fix $a>1$ and let $\nu$ be the probability measure on $\{1,2,\ldots\}$ defined by $\nu_k=k^{-a}/\zeta(a)$, where $\zeta(a)=\sum_{k=1}^\infty k^{-a}$ denotes the Riemann zeta function. Since relative entropy  $\sum_k \mu_k \log (\mu_k/\nu_k)$ between any two probability measures is non-negative, we have that
\[
 \sum_{k=1}^\infty \mu_k \log \frac{\mu_k}{\nu_k} = -\operatorname{Ent}(\mu)+a\sum_{k=1}^\infty \mu_k \log k + \log \zeta(a) \geq 0.
\]
Rearranging this inequality yields that
\[
  \sum_{k=1}^\infty \mu_k \log k \geq \frac{1}{a}\operatorname{Ent}(\mu)-\frac{1}{a}\log \zeta(a),
\]
which is easily seen to imply the claim since $\frac{1}{2}\log \zeta(2)=0.24885\ldots \leq 1/4$.
\end{proof}

We now apply \cref{lem:entropy_alternative} to prove \cref{prop:entropy_and_cell_size}.

\begin{proof}[Proof of \cref{prop:entropy_and_cell_size}] We have by \cref{lem:coupling_two_point} that
\[
\frac{1}{2}\sum_{y\in S} \left[\frac{\kappa(y,z)}{\kappa(x,z)} \wedge 1\right] \le
\E\left[\#[x] \mid X(x)=z\right] \leq  \sum_{y\in S} \left[\frac{\kappa(y,z)}{\kappa(x,z)} \wedge 1\right]\]
 and hence by the definition of $C$-reversibility that
\[
\frac{1}{2C^2} \sum_{y\in S} \left[\frac{\kappa(z,y)}{\kappa(z,x)} \wedge 1 \right] \leq \E\left[\#[x] \mid X(x)=z\right] \leq C^2 \sum_{y\in S} \left[\frac{\kappa(z,y)}{\kappa(z,x)} \wedge 1 \right]
\]
for each $x,z\in S$.
Using Jensen's inequality and $C$-reversibility a second time, we deduce that
\begin{align*}
  \E \log \#[x] \leq  \E \left[\log \E\left[\#[x] \mid X(x)\right]\right] &\leq 
  \sum_{z} \kappa(x,z) \log C^2 \sum_{y\in S} \left[\frac{\kappa(z,y)}{\kappa(z,x)} \wedge 1 \right]\\
  &\leq  C
  \sum_{z} \kappa(z,x) \log \sum_{y\in S} \left[\frac{\kappa(z,y)}{\kappa(z,x)} \wedge 1 \right] + 2\log C
\end{align*}
and hence by \cref{lem:entropy_alternative} that
\[
  \frac{1}{|S|}\sum_{x\in S} \E \log \#[x] \leq \frac{C}{|S|}\sum_{x,z\in S} \kappa(z,x) \log \sum_{y\in S} \left[\frac{\kappa(z,y)}{\kappa(z,x)} \wedge 1 \right] + 2\log C \leq  \frac{C}{|S|}\sum_{z\in S} \operatorname{Ent}(\mu_z)+2\log C.
\]
It remains to prove a matching lower bound.
On the event $X(x)=z$, each $y$ with $\kappa(z,y) \geq \kappa(z,x)$ satisfies 
\[
  \kappa(y,z) \geq C^{-1} \kappa(z,y) \geq C^{-1} \kappa(z,x) \geq C^{-2} \kappa(x,z)
\]
and therefore has $X(y)=z$ with probability at least $\frac{1}{2}C^{-2}$ by \cref{lem:coupling_two_point}.
 Applying
Markov's inequality\footnote{More precisely we are using the fact that if $X$ is a random variable taking values in $[0,M]$ with $\E X \geq \theta M$ then $\P(X \leq \frac{1}{2}\theta M) = \P(M-X \geq (1-\frac{1}{2}\theta)M ) \leq (1-\theta)/(1-\frac{1}{2}\theta) \leq 1-\frac{1}{2}\theta$.}
 to the conditional law of $\#\{y: \kappa(z,y) \geq \kappa(z,x), y\notin [x] \}$ given $X(x)=z$, we obtain that
\[
  \P\left(\#[x] \geq \frac{1}{4C^2}\#\{y: \kappa(z,y) \geq \kappa(z,x) \} \;\middle|\; X(x)=z\right) \geq \frac{1}{4C^2}
\]
and hence that
\begin{align*}
\E \left[\log \#[x] \mid X(x)=z\right] \geq \frac{1}{4C^2} \log \frac{1}{4C^2}\#\{y: \kappa(z,y) \geq \kappa(z,x) \}.
\end{align*}
Taking expectations over $X(x)$ and averaging over the choice of $x\in S$ we obtain that
\begin{align*}
  \frac{1}{|S|}\sum_{x\in S}\E \log \# [x] &\geq \frac{1}{4C^2} \frac{1}{|S|}\sum_{x,z\in S}\kappa(x,z) \left[\log \#\{y: \kappa(z,y) \geq \kappa(z,x)\}- \log 4C^2\right]\\
  &\geq \frac{1}{4C^3}\frac{1}{|S|}\sum_{x,z\in S}\kappa(z,x) \log \#\{y: \kappa(z,y) \geq \kappa(z,x)\}-\frac{\log 4C^2}{4C^2}
\\&\geq  \frac{1}{8C^3}\frac{1}{|S|}\sum_{z\in S} \operatorname{Ent}(\mu_z)-\frac{1}{16C^3}-\frac{\log 4C^2}{4C^2},
\end{align*}
where we used $C$-reversibility again in the second line and used \cref{lem:entropy_alternative} in the last line.
The claimed lower bound follows since the additive constant appearing here is at most $(1+8\log 2)/16 \leq 1/2$ for $C\geq 1$.
\end{proof}

We now state and prove an analogue of \cref{prop:entropy_and_cell_size} for infinite unimodular transitive graphs.

\begin{prop}[Entropy is equivalent to cell size; transitive case]
\label{prop:entropy_and_cell_size_transitive}
Let $G=(V,E)$ be a locally finite unimodular transitive graph, let $\kappa:V\times V\to [0,1]$ be an automorphism-invariant probability kernel on $G$, and consider the grand Poisson cell decomposition associated to $\kappa$. If $\kappa$ is $C$-reversible  for some $C\geq 1$ and we define $\mu_x(\,\cdot\,)=\kappa(x,\,\cdot\,)$ for each $x\in V$ then 
\[
 \frac{1}{8C^3} \operatorname{Ent}(\mu_o) -\frac{1}{2} \leq \E\log \#[o] \leq C \operatorname{Ent}(\mu_o)+2\log C.
\] 
\end{prop}

\begin{proof}[Proof of \cref{prop:entropy_and_cell_size_transitive}]
We have as in the proof of \cref{prop:entropy_and_cell_size}  that
\[
  \E \log \# [o] \leq C \sum_{z\in V} \kappa(z,o) \log \sum_{y\in V} \left[\frac{\kappa(z,y)}{\kappa(z,o)} \wedge 1  \right] + 2\log C
\]
and
\[
  \E \log \#[o] \geq \frac{1}{4C^3} \sum_{z\in V} \kappa(z,o) \log  \# \{y \in V : \kappa(z,y) \geq \kappa(z,o)\} - \frac{\log 4C^2}{4C^2}.
\]
Indeed, both inequalities hold for any $C$-reversible probability kernel on any countable set $V$ and any $o\in V$. Applying the mass-transport principle 
with the functions
\[
  (u,v) \mapsto \kappa(v,u) \log \sum_{y\in V} \left[\frac{\kappa(v,y)}{\kappa(v,u)} \wedge 1  \right] \quad \text{and} \quad (u,v)\mapsto \kappa(v,u) \log \# \{y \in V : \kappa(v,y) \geq \kappa(v,u)\},
\]
it follows from this and \cref{lem:entropy_alternative} that
\[
  \E \log \# [o] \leq C \sum_{z\in V} \kappa(o,z) \log \sum_{y\in V} \left[\frac{\kappa(o,y)}{\kappa(o,z)} \wedge 1  \right] + 2\log C \leq C \operatorname{Ent}(\mu_o)+2\log C
\]
and
\begin{multline*}
  \E \log \#[o] \geq \frac{1}{4C^3} \sum_{z\in V} \kappa(o,z) \log  \# \{y \in V : \kappa(o,y) \geq \kappa(o,z)\} - \frac{\log 4C^2}{4C^2}
  \\
  \geq \frac{1}{8C^3} \operatorname{Ent}(\mu_o) - \frac{1}{16C^3} - \frac{\log 4C^2}{4C^2} \geq \frac{1}{8C^3} \operatorname{Ent}(\mu_o) - \frac{1}{2}
\end{multline*}
as claimed.
\end{proof}

\section{Proof of  the main theorems for transitive graphs}\label{sec:proof_of_the_main_theorem}

In this section we prove our main theorem concerning transitive graphs, \cref{thm:growth_main}.
%
%
Let $G=(V,E)$ be a transitive graph of non-negative Ollivier--Ricci curvature which we fix for the remainder of this section. We let $o\in V$ be a fixed root vertex and write $\operatorname{Gr}(r)=\# B(o,r)$ for the volume of the ball of radius $r$ around $o$. We also write 
$D_n = \E_o d(o,X_n)$ for the expected $n$-step displacement of the random walk started at $o$.
We wish to prove that
\[
  \log \operatorname{Gr}(r) \preceq d^3 (\log r),
\]
where we recall that we write $\asymp$, $\preceq$, and $\succeq$ for equalities and inequalities holding to within universal positive multiplicative constants.

\medskip

For each $n\geq 1$ the $n$-step lazy random walk transition matrix $P^n:V^2\to[0,1]$ is a $1$-reversible automorphism-invariant probability kernel in the sense of \cref{sec:the_grand_poisson_cell_decomposition}. The proof of \cref{prop:log_squared_transitive} will proceed by studying the grand Poisson cell decomposition associated to a certain truncated family of random walk measures obtained by conditioning the random walk to have displacement not too much larger than $D_n$.
We will need to be careful precisely how this truncation is performed, with the existence of a ``good scale for truncation'' guaranteed by the following elementary lemma.


\begin{lemma}
\label{lem:truncation_scale}
For each $n\geq 1$ there exists $\Lambda_n \in [1,10]$ such that
\[
  \mathbf{P}_x\Bigl(X_n \in B(x,\lceil 2D_n+\Lambda_nn^{1/2}\rceil) \setminus B(x,\lceil 2D_n+\Lambda_nn^{1/2}\rceil -2)\Bigr) \leq \frac{1}{\sqrt{n}}
\]
for every $x\in V$ and $n\geq 1$.
\end{lemma}

\begin{proof}[Proof of \cref{lem:truncation_scale}]
This is a trivial consequence of the fact that
\[
\sum_{k=1}^\infty \mathbf{P}_x\Bigl(X_n \in B(x,2k) \setminus B(x,2k -2)\Bigr) \leq 1\]
by additivity of probability measures.
\end{proof}

For each $n\geq 1$ we fix $\Lambda_n \in [1,10]$ as guaranteed to exist by \cref{lem:truncation_scale} and define the probability kernel $\widehat P^n$ by conditioning the lazy random walk to have displacement at most $\lceil 2D_n + \Lambda_n n^{1/2}\rceil$. The probability kernel $\widehat P^n$ is $1$-reversible and automorphism-invariant. Note also that the random walk on $G$ satisfies
\[
  \mathbf{P}_o(d(o,X_n) \geq \lceil 2D_n +\Lambda_n n^{1/2} \rceil )  \leq \frac{1}{2}
\]
by Markov's inequality, so that if $\mu_x^n$ and $\widehat{\mu}_x^n$ denote the probability measures with mass functions $P^n(x,\cdot)$ and $\widehat P^n(x,\cdot)$ then $\widehat{\mu}_x^n(y)\leq 2 \mu_x^n(y)$ for every $x,y\in V$ and $n\geq 1$.

\medskip

We next prove that these truncated measures satisfy a total variation bound just as strong as that established for the untruncated measures $(\mu_x^n)_{x\in V}$ in \cref{thm:Salez}.

\begin{lemma} The estimate
\label{lem:mu_hat_TV}
\begin{equation}\notag
      \|\widehat\mu_x^n-\widehat\mu_y^n\|_{\operatorname{TV}}\preceq{\sqrt{\frac{d}{n}}d(x,y)}
\end{equation}
holds for every $x,y \in V$ and $n\geq 1$.
\end{lemma}

\begin{proof}[Proof of \cref{lem:mu_hat_TV}]
It suffices to prove the inequality in the case $x\sim y$. Writing
$R= \lceil 2D_n + \Lambda_n n^{1/2}\rceil$ and
 $B=B(x,R) \cap B(y,R) $, we can write
\begin{multline*}
  \|\widehat\mu_x^n-\widehat\mu_y^n\|_{\operatorname{TV}} = \frac{1}{2}\sum_{z\in B} |\widehat\mu_x^n(z)-\widehat\mu_y^n(z)| + \frac{1}{2}\sum_{z\in V \setminus B} |\widehat\mu_x^n(z)-\widehat\mu_y^n(z)| 
  \\
  \leq 2 \|\mu_x^n-\mu_y^n\|_\mathrm{TV} +\frac{1}{2} \widehat\mu_x^n(V\setminus B) + \frac{1}{2}\widehat\mu_y^n(V\setminus B),
\end{multline*}
where we used that if $z\in B$ then $\widehat\mu_x^n(z)-\widehat\mu_y^n(z)$ is proportional to $\mu_x^n(z)-\mu_y^n(z)$ with constant of proportionality between $1$ and $2$.
Noting that $B$ contains the ball of radius $R-1=\lceil 2D_n + \Lambda_n n^{1/2}\rceil-1$ around $x$, it follows by the choice of $\Lambda_n$ that $\widehat\mu_x^n(V\setminus B),\widehat\mu_y^n(V\setminus B) \leq 2 n^{-1/2}$, and the claim follows from this and \cref{thm:Salez}.
\end{proof}

Now we define $\widehat{[x]}_n$ to be the cell corresponding to $x$ in the grand Poisson decomposition with respect to $(\widehat{\mu}_x^n)_{x}$, where we take the Poisson point processes defining $\widehat{[x]}_n$ and $\widehat{[x]}_m$ to be independent for distinct values of $n$ and $m$. We have by construction that
\begin{equation}\label{eq:hat_cell_diameter}
  \diam(\widehat{[o]}_n)\leq 4D_n+20\sqrt{n}+2
\end{equation}
for every $n\geq 1$.
We define 
\begin{equation}\notag
  \widehat{a}_n    \defeq \E \log \# {\widehat{[o]}}_n
\end{equation}
and seek to establish functional inequalities relating $\widehat{a}_n$, $D_n$, and $\log \operatorname{Gr}(r)$.
The simplest such inequality arises as a direct consequence of \cref{lem:mu_hat_TV} and the defining properties of the grand Poisson cell decomposition.




\begin{lemma}\label{lem:Gr_bounded_above_by_hat_a}
  There exists a universal  constant $ A$ such that
  \begin{equation}\notag
    \log\Gr(r)\leq 2\widehat{a}_{n}+\log 2
  \end{equation}
  for every $r\geq 1$ and $n\geq Adr^2$.
\end{lemma} 

\begin{proof}[Proof of \cref{lem:Gr_bounded_above_by_hat_a}]
The property \eqref{eq:coupling_TV} of the grand Poisson cell decomposition and the total variation estimate  of \cref{lem:mu_hat_TV} imply that if $x\in  B(o,r)$ then 
\begin{equation}\notag
\P(x\notin\widehat{[o]}_{n})\leq 2    \|{\widehat\mu_x^{n}}-{\widehat\mu_o^{n}}\|_\TV \preceq  \sqrt{\frac{d}{n}}r .
\end{equation}
As such, there exists a universal constant $A$ such that if $n,r\geq 1$ are such that $n\geq Adr^2$ then 
$\P(x\notin{\widehat{[o]}}_{n})\leq \frac{1}{4}$ for every $x\in B(o,r)$. Summing over $x\in B(o,r)$ yields that
\begin{equation*}
  \E\Big(\#\big(B(o,r)\setminus \widehat{[o]}_{n}\big)\Big)\leq \frac{1}{4}\# B(o,r),
\end{equation*}
whence, by Markov's inequality,
\begin{equation}\notag
  \P\Big(\#\big(B(o,r)\cap \widehat{[o]}_{n}\big)\leq \frac{1}{2}\# B(o,r)\Big)=  \P\Big(\#\big( B(o,r)\setminus \widehat{[o]}_{n}\big)> \frac{1}{2}\# B(o,r)\Big)\leq \frac{1}{2}.
\end{equation}
Therefore, 
\begin{equation}\notag
  \E \log \#  \widehat{[o]}_n\geq \frac{1}{2}\log\frac{\# B(o,r)}{2},
\end{equation}
which is exactly the claim.
\end{proof}

Next, we use the isoperimetric inequality \cref{thm:munch} to derive the following relationship between $\widehat{a}_n$ and $D_n$.

\begin{prop}[Recurrence relation for $   \widehat a_n $]\label{prop:hat_a_functional_inequality_with_D_n}
It holds that 
  \begin{equation}\notag
    \widehat a_n-\widehat a_m \preceq   \frac{d^{3/2}}{\sqrt{m}}(D_n+\sqrt{n}) +1
  \end{equation}
  for every $n\geq  m \geq 1$.
\end{prop}

\begin{proof}[Proof of \cref{prop:hat_a_functional_inequality_with_D_n}]
  We apply the isoperimetric  inequality of \Cref{thm:munch} in the form of \eqref{sdcdsacsc} to the cells $\widehat{[o]}_{m,n}\defeq {\widehat{[o]}}_m\cap \widehat{[o]}_n\subseteq \widehat{[o]}_n$ to obtain that
  \begin{equation}\notag
    \log\#\widehat{[o]}_n\leq \log\#\widehat{[o]}_{m,n}+4 \diam(\widehat{[o]}_n)\frac{\#\partial\widehat{[o]}_{m,n}}{\#\widehat{[o]}_{m,n}}+\log 2
  \end{equation}
  for every $n,m \geq 1$. 
   Bounding $\#\widehat{[o]}_{m,n} \leq \#\widehat{[o]}_{m}$, taking expectations and using the bound \eqref{eq:hat_cell_diameter}, 
    we obtain that
\begin{equation}\notag
         \widehat a_n- \widehat a_m\preceq   (D_n+\sqrt{n})\E\bigg[ \frac{\#\partial\widehat{[o]}_{m,n}}{\#\widehat{[o]}_{m,n}}\bigg]+1.
  \end{equation}
Applying \Cref{lem:cell_isoperimetry_transitive}, the boundary-to-volume ratio appearing here can be bounded
  \begin{align*}
      \E\bigg[ \frac{\#\partial\widehat{[o]}_{m,n}}{\#\widehat{[o]}_{m,n}}\bigg]&=\E \sum_{z\sim o}\mathbbm{1}(\widehat{[z]}_{m,n}\ne \widehat{[o]}_{m,n})\leq \E \sum_{z\sim o}\mathbbm{1}(\{\widehat{[z]}_{m}\ne \widehat{[o]}_{m}\})+
      \E \sum_{z\sim o}\mathbbm{1}(\{\widehat{[z]}_{n}\ne \widehat{[o]}_{n}\})\\&\leq \sum_{x\sim o} \frac{2\|{\widehat{\mu}_x^m}-{\widehat{\mu}_o^m}\|_\TV}{1+\|{\widehat{\mu}_x^m}-{\widehat{\mu}_o^m}\|_\TV}+\sum_{x\sim o} \frac{2\|\widehat{\mu}_x^n-\widehat{\mu}_o^n\|_\TV}{1+\|\widehat{\mu}_x^n-\widehat{\mu}_o^n\|_\TV}
      \preceq 
      d\sqrt{\frac{d}{m}}
  \end{align*} 
 for every $n\geq m \geq 1$ as required, where we applied \cref{lem:mu_hat_TV} in the last inequality.
\end{proof}

Our strategy for the remainder of the section is to find relations between $\log \operatorname{Gr}(r)$, $D_n$, and $\widehat{a}_n$ making \cref{lem:Gr_bounded_above_by_hat_a,prop:hat_a_functional_inequality_with_D_n} into a functional inequality strong enough to force polynomial growth. We do this in two stages, first using the Varopoulos--Carne inequality to prove a functional inequality strong enough to force the quasi-polynomial estimate $\log \operatorname{Gr}(r)=O((\log r)^2)$ then using \cref{prop:displacement_doubling} to prove a different functional inequality that, given this quasi-polynomial bound as an input, is strong enough to force polynomial growth.

\medskip

\noindent \textbf{A quasi-polynomial upper bound.}
Our first functional inequality is derived from \cref{prop:hat_a_functional_inequality_with_D_n} together with the following relationship between the expected displacement and $\widehat{a}_n$ holding for any transitive graph. 
This relationship is a manifestation of the well-known inequality $D_n^2 \preceq n H_n$ between the expected displacement and the Shannon entropy (closely related to the so-called ``fundamental inequality'' for random walks on groups) together with the relationship $\E \log \#[o]_n \asymp H_n$ provided by \cref{prop:entropy_and_cell_size_transitive}.


\begin{lemma}
\label{lem:D_n_upper}
The inequality
    \begin{equation}\notag
  D_n^2\preceq n \widehat{a}_n +n\log n +\log (d+1)
\end{equation}  
holds for every $n\geq 1$.
\end{lemma}

The proof of this lemma will rely on the  Varopoulos--Carne inequality \cite{Varopoulos1985b,Carne1985} (see also \cite[Theorem 13.4]{LP:book}), which states that the lazy random walk on any graph satisfies the Gaussian-like pointwise estimate
 \[
  P^n(x,y) \leq 2 \sqrt{\frac{\deg(y)}{\deg(x)}}  \exp\Big(-\frac{d(x,y)^2}{2n}\Big)\qquad\text{for every $x,y\in V$ and $n\geq 1$}.
  \]
For transitive graphs this inequality allows us to lower bound the Shannon entropy of the random walk in terms of its average displacement, yielding that
  \begin{multline}
  \label{eq:displacement_v_entropy}
   H_n\defeq -\sum_{y\in V}P^n(o,y)\log P^n(o,y)\geq \sum_{y\in V}P^n(o,y)\frac{d(o,y)^2}{2n}-\log 2
    \\=\frac{1}{2n}\mathbf{E}_o d(o,X_n)^2-\log 2 \geq \frac{1}{2n}D_n^2-\log 2
  \end{multline}
where we applied Jensen's inequality to bound the mean-squared displacement by the squared mean displacement $D_n^2$. We will also use the fact that
if $\mu$ is a probability measure on a countable set $S$ and $\mu|_A$ is defined from $\mu$ by conditioning on some set $A$ with $\mu(A)>0$ then
\begin{equation}\notag
  \|\mu-\mu|_A\|_\mathrm{TV} = \frac{1}{2}\sum_{x\in A} \left[\frac{\mu(x)}{\mu(A)}-\mu(x)\right]  + \frac{1}{2}\sum_{x\in S\setminus A} \mu(x) =  \mu(S\setminus A).
\end{equation}


\begin{proof}[Proof of \cref{lem:D_n_upper}]
We may assume that $D_n \geq 4\sqrt{n\log n}$, the claimed inequality holding trivially otherwise. Let $[o]_n$ and $\widehat{[o]}_n$ be the cells of the origin in the grand Poisson cell decompositions associated to $P^n$ and $\widehat{P}^n$ computed with respect to the same underlying Poisson process.
It follows from \eqref{eq:displacement_v_entropy} and \cref{prop:entropy_and_cell_size_transitive} that
    \begin{equation}\label{tc1}
    D_n^2\preceq n (\operatorname{Ent}(\mu^n_o)+\log 2) \preceq  n (\E \log \#[o]_n + 1)
  \end{equation}
  where $[o]_n$ denotes the cell of $o$ associated to the grand Poisson decomposition with respect to the probability kernel $P^n$.
  On the other hand, we have by Markov's inequality, \Cref{lem:different_cells},  and \cref{prop:Lipschitz_concentration} applied with $f=d(o,\,\cdot\,)$ that
  \begin{multline*}
    \P \bigg( \#\widehat{[o]}_n\leq \frac{1}{2}\#{[o]_n}\bigg)\leq 
    \P \bigg( \#\big({[o]}_n\setminus \widehat{[o]}_n\big)\geq \frac{1}{2}\#{[o]_n}\bigg)
    \leq 2\E\bigg(\frac{\#\big([o]_n\setminus\widehat{[o]}_n\big)}{\#{[o]}_n}\bigg)
    \\\leq 8\|\mu_o^n-\widehat\mu_o^n\|_\TV
    \leq 8 \mathbf{P}_o(d(o,X_n) \geq 2D_n + \Lambda_n n^{1/2}) 
    \leq 16 \exp\left[- \frac{D_n^2}{8n}\right] \leq 16 n^{-2}
  \end{multline*}
  where we used the assumption that $D_n \geq 4\sqrt{n\log n}$ in the final inequality. Since $\diam ([o]_n)\leq 2n$, it follows that
  \begin{multline}\label{tc3}
      \E \log\#[o]_n\leq \E\bigg[ \log\#[o]_n\mathbbm{1}(\#\widehat{[o]}_n\geq \frac{1}{2}\#{[o]_n})\bigg]+  \E \bigg[\log\#[o]_n\mathbbm{1}( \#\widehat{[o]}_n\leq \frac{1}{2}\#{[o]_n})\bigg]\\
    \preceq \E \log (2\#\widehat{[o]}_n) + \log \Gr(2n) n^{-2}
    \preceq \E \log \#\widehat{[o]}_n+\frac{\log (d+1)}{n}+ 1,
  \end{multline}
  where we used the trivial bound $\operatorname{Gr}(r) \leq \sum_{i=0}^r d^i \preceq (d+1)^{r+1}$ in the final inequality.
  We conclude the proof by combining \eqref{tc1} and \eqref{tc3}. (We have an $n\log n$ rather than an $n$ in the statement so that the inequality also holds if our initial assumption that $D_n \geq 4 \sqrt{n\log n}$ is false.)
\end{proof}

We now analyze the functional inequality for $\widehat{a}_n$ obtained from \cref{prop:hat_a_functional_inequality_with_D_n,lem:D_n_upper} to obtain the following quasi-polynomial growth bound.

\begin{prop}\label{prop:log_squared_transitive}It holds that
$\log \Gr(r)\preceq d^3 (\log r)^2$
  for every $r\geq d+1$.
\end{prop}

  \begin{proof}[Proof of \cref{prop:log_squared_transitive}]
    We insert the inequality of \cref{lem:D_n_upper} into that of \cref{prop:hat_a_functional_inequality_with_D_n}  to obtain that 
    \begin{equation}\notag
             \widehat a_n-   \widehat a_m \preceq \big( \sqrt{n\widehat{a}_n+n\log n +\log (d+1)}+\sqrt{n}\big)\frac{d^{3/2}}{\sqrt{m}} +1
    \end{equation}
    for every $n\geq m$
and hence that
    \begin{equation}\label{eq:functional_inequality01}
             \widehat a_n-   \widehat a_m \preceq \sqrt{\frac{d^3 n\widehat a_n}{m}} 
    \end{equation}
    for every $n\geq 2m \geq 2d$ such that $\widehat{a}_n \geq \log n$.
    (That is, we note that if $\widehat{a}_n \geq \log n$ then the various nuisance terms can be absorbed into the main term on the right at the cost of increasing the implicit constant.)
    Defining the sequence $(b_k)_{k=1}^\infty$ as $b_k \defeq  \widehat a_{d 2^k}^{1/2}$, 
    it follows in particular that there exists a universal constant $A_1$ such that if $b_k \geq \max\{k, A_1 d^{3/2}\}$ for some $k\geq 2$ (which guarantees $\widehat{a}_{d2^k} \geq \log (d2^k)$) then
    $b_k \leq 2 b_{k-1}$
and hence by a second application of \eqref{eq:functional_inequality01} that
\[
  b_k^2 - b_{k-1}^2  \preceq  d^{3/2} b_{k-1}.
\]
Using that $ b_k^2 - b_{k-1}^2 = (b_k-b_{k-1})(b_k+b_{k-1}) \geq (b_k-b_{k-1}) b_{k-1}$ when $b_k\geq b_{k-1}$, we obtain that if $b_k \geq \max\{k, A_1 d^{3/2}\}$ then
\[
  b_k -b_{k-1} \preceq d^{3/2},
\]
this inequality holding trivially if $b_k<b_{k-1}$.
As such, if for each $k\geq 1$ we define $k_0 = \max\{ 1\leq \ell \leq k : \ell=1$ or $b_\ell < \max\{\ell, A_1 d^{3/2}\}\}$ then
\[
  b_k \leq b_{k_0}+\sum_{\ell=k_0+1}^{k} (b_\ell -b_{\ell-1}) \preceq k d^{3/2},
\]
for every $k\geq 1$, 
where we used that
 every cell of the decomposition associated to $\widehat P^{2d}$ has diameter at most $4d$ and hence that
 $b_1 =\widehat{a}_{2d}^{1/2} \leq (\log \operatorname{Gr}(4d))^{1/2} \preceq \sqrt{d \log (d+1)} \preceq d^{3/2}$ to handle the case $k_0=1$.
We deduce that $\widehat{a}_{d2^k} \preceq d^3 k^2$ for every $k\geq 1$, and the claimed upper bound on $\log \operatorname{Gr}(r)$ follows immediately from this and \cref{lem:Gr_bounded_above_by_hat_a}.
\end{proof}

\noindent \textbf{Improving to a polynomial bound.}
Our next goal is to improve the quasi-polynomial bound of \cref{prop:log_squared_transitive} to a polynomial bound. The main technical step in this part of the proof is the following functional inequality for the growth derived from \cref{prop:hat_a_functional_inequality_with_D_n} and the displacement bound \cref{prop:displacement_doubling}.

\begin{prop}[Functional inequality for $\log\Gr$]\label{prop:main_transitive_functional_inequality} 
The inequality
  \begin{equation}\label{eq:main_transitive_functional_inequality}
  \log \Gr(r)\preceq d^{3/2}\sqrt{\log r\log\Gr(A dr^2)}+ d^{3/2}\log r
\end{equation}  
holds for every $r\geq d+1$.
\end{prop}

Note that this inequality is not useful without the weaker bound on the growth provided by \cref{prop:log_squared_transitive} (indeed, the inequality holds vacuously for graphs of exponential volume growth).

\begin{proof}[Proof of \cref{prop:main_transitive_functional_inequality}]
  For each $k\geq 1$ we set
  \begin{equation}\notag
    m_k\defeq\max\{ 0\leq m\leq k :D_{2^k}/2\leq D_{2^{k-m}} \leq 2 D_{2^k} \}.
  \end{equation}
  We have from \Cref{prop:displacement_doubling} that for every $k\geq 1$,
  \begin{equation}\label{dscascsc}
      D_{2^k} \leq 2D_{2^{k-{m_k}}} \preceq \sqrt{2^{k-m_k} \log \frac{5 {\Gr}(4D_{2^{k-m_k}})}{{\Gr}(D_{2^{k-m_k}}/8)}}\preceq \sqrt{2^{k-m_k} \log \frac{5 {\Gr}(8D_{2^{k}})}{{\Gr}(D_{2^{k}}/16)}}.
  \end{equation}
This estimate will allow us to recognize a telescoping series in the functional inequality, with the extra $2^{-m_k/2}$ letting us handle the possibility that the sequence $k\mapsto D_{2^k}$ takes values of the same order many times.
More specifically,
we will use the fact that if we define $L_i\defeq \{\ell \geq 0:D_{2^\ell}\in [2^{i-1},2^{i}] \}$ for each $i\geq 1$ then
\begin{equation}
\label{eq:L_i_dyadic}
  \sum_{\ell \in L_i} 2^{-m_\ell} \leq 2
\end{equation}
for every $i\geq 1$. This inequality holds vacuously if $L_i = \emptyset$. Otherwise, if we let $\ell_0(i)$ be the smallest element of $L_i$ and note that if $\ell\in L_i$ then $m_\ell \geq \ell-\ell_0(i)$  so that
  $\sum_{\ell \in L_i} 2^{-m_\ell} \leq \sum_{\ell \geq \ell_0(i)} 2^{-(\ell-\ell_0(i))} \leq 2$
as claimed.

\medskip

The inequality in  \Cref{prop:hat_a_functional_inequality_with_D_n} (applied with $n=2^k$ and $m={2^{k-1}}=n/2$) yields that
\begin{align*}
 \widehat{a}_{2^k}-\widehat{a}_{2^{k-1}}\preceq d^{3/2}\frac{D_{2^k}}{2^{k/2}}+ d^{3/2}
\end{align*}
for every $k\geq 2$.
Since  $\diam(\widehat{[o]}_2)\leq 4$ and hence $\widehat{a}_2\preceq \log(d+1) \preceq d^{3/2}$, we can sum this estimate over $k$ to obtain that
\begin{equation}\notag
    \widehat{a}_{2^k}
    \preceq  d^{3/2}\sum_{\ell=1}^{k}\frac{D_{2^\ell}}{2^{\ell/2}}+d^{3/2}k    \preceq d^{3/2}\sum_{\ell=1}^{k}\mathbbm{1}(D_{2^\ell}\geq 1)\frac{D_{2^\ell}}{2^{\ell/2}}+d^{3/2}k.
\end{equation}
To estimate the sum in the inequality above, we use \eqref{dscascsc} and Cauchy--Schwarz inequality to obtain that
\begin{multline*}
    \sum_{\ell=1}^{k}\mathbbm{1}(D_{2^\ell}\geq 1)\frac{D_{2^\ell}}{2^{\ell/2}}\preceq \sum_{\ell=1}^{k} \sqrt{\mathbbm{1}(D_{2^\ell}\geq 1) 2^{-m_\ell}\log \frac{5 {\Gr}(8D_{2^{\ell}})}{{\Gr}(D_{2^{\ell}}/16)}}
    \\\preceq  \sqrt{ k\sum_{\ell=1}^{k}\mathbbm{1}(D_{2^\ell}\geq 1)2^{-m_\ell}\log \frac{5 {\Gr}(8D_{2^{\ell}})}{{\Gr}(D_{2^{\ell}}/16)}}
\preceq k + \sqrt{ k\sum_{\ell=1}^{k}\mathbbm{1}(D_{2^\ell}\geq 1)2^{-m_\ell}\log \frac{{\Gr}(8D_{2^{\ell}})}{{\Gr}(D_{2^{\ell}}/16)}}
\end{multline*}
and hence that
\begin{equation}\label{cdscsacdc}
    \widehat{a}_{2^k}\preceq d^{3/2}\sqrt{ k\sum_{\ell=1}^{k}\mathbbm{1}(D_{2^\ell}\geq 1)2^{-m_\ell}\log \frac{{\Gr}(8D_{2^{\ell}})}{{\Gr}(D_{2^{\ell}}/16)}}+d^{3/2}k
\end{equation}
for every $k\geq 1$.
We now use the inequality \eqref{eq:L_i_dyadic} to bound the sum appearing in  \eqref{cdscsacdc}. Noting that $\log_2 D_{2^\ell} \leq k$ for every $\ell\leq k$, we have that
\begin{align*}
  \sum_{\ell=1}^{k} \mathbbm{1}(D_{2^\ell}\geq 1) 2^{-m_\ell} \log \frac{{\Gr}(8D_{2^\ell})}{{\Gr}(D_{2^\ell}/16)}
  &\leq \sum_{\ell=1}^{k} \sum_{i=1}^{k}\mathbbm{1}(D_{2^\ell}\in [2^{i-1},2^{i}])2^{-m_\ell}\log \frac{{\Gr}(8D_{2^\ell})}{{\Gr}(D_{2^\ell}/16)}\\
  &\leq   \sum_{i=1}^{k}\sum_{\ell=1}^{k} \mathbbm{1}(\ell\in L_{i})2^{-m_\ell}\log \frac{{\Gr}(2^{i+3})}{{\Gr}(2^{i-5})}
\end{align*}
and hence by \eqref{eq:L_i_dyadic} that
\begin{multline*}
  \sum_{\ell=1}^{k} \mathbbm{1}(D_{2^\ell}\geq 1) 2^{-m_\ell} \log \frac{{\Gr}(8D_{2^\ell})}{{\Gr}(D_{2^\ell}/16)}
  \leq 2  \sum_{i=1}^{k}\log \frac{{\Gr}(2^{i+3})}{{\Gr}(2^{i-5})}\\
  = 2\sum_{i=1}^{k}\log {{\Gr}(2^{i+3})}-  2\sum_{i=1}^{k}\log{{\Gr}(2^{i-5})}
  \preceq \log\Gr(2^{k+4}).
\end{multline*}
Substituting this inequality into \eqref{cdscsacdc}, we obtain that
\begin{equation}\label{cdsxcasdc}
      \widehat{a}_{2^k}\preceq d^{3/2} \sqrt{k\log\Gr(2^{k+4})}+d^{3/2}k 
\end{equation}
for every $k\geq 1$, where we used the running assumption that our graph is not an isolated vertex to bound $\log\Gr(2^{k+4}) \succeq 1$.

\medskip

Now, by \cref{lem:Gr_bounded_above_by_hat_a}, there exists a universal constant $A$ such that 
\begin{equation}\notag
  \log\Gr(r)\leq 2\widehat{a}_{n}+\log 2\qquad\text{for every $r\geq 1$}
\end{equation}
for every $n\geq A d r^2$.
Taking $n$ to be an integer power of $2$ between $Adr^2$ and $2Adr^2$, it follows from \eqref{cdsxcasdc} that
%
%
\begin{align*}
  \log \Gr(r)&\preceq d^{3/2} \sqrt{\log(8 A dr^2)\log\Gr(2^5Adr^2) }+ d^{3/2}\log(8 Adr^2)
  \end{align*}
for every $r\geq 1$ and hence that
\begin{align*}
  \log \Gr(r)&\preceq d^{3/2} \sqrt{\log(r)\log\Gr(2^5Adr^2) }+ d^{3/2}\log(r)
  \end{align*}
for every $r \geq d+1$, which easily implies the claim.
%
\end{proof}

\begin{proof}[Proof of \cref{thm:growth_main}]
It suffices to prove the claim for $r\geq d+1$, the inequality holding trivially for smaller values of $r\geq 2$ since $\operatorname{Gr}(r)\leq \sum_{i=0}^r d^i \preceq (d+1)^{r+1}$.
By \eqref{eq:main_transitive_functional_inequality}, there exist universal constants $C,A \geq 1$ such that
    \begin{equation}\label{dscsacsac}
    \log \Gr(r)\leq  C d^{3/2}\big(\sqrt{\log r\log\Gr(A d r^2)}+\log r\big)
  \end{equation}  
  for every $r\geq d+1$. On the other hand, we have by \cref{prop:log_squared_transitive} that there exists a universal constant $B\geq 1$ such that
  \begin{equation}
    \log \operatorname{Gr}(r) \leq B d^3 (\log r)^2
  \label{eq:weak_growth_restate}
  \end{equation}
  for every $r\geq d+1$.   Let $B'$ be the universal constant defined by
  \[
   B' = \max \left\{B, 4C^2 \max_{a \geq b+1 \geq 2} \left(\frac{\log(A ba^2)}{\log a}\right)^2 \right\},
   \]
   so that $B'\geq B$ and, for $\alpha\in [0,1]$,
\[
    2C \sqrt{B' \log r (\log (A dr^2))^{1+\alpha}} \leq B' (\log r)^{1+\alpha/2}
  \]
  for every $r \geq d+1 \geq 2$.
We claim that
\[
  \log \operatorname{Gr}(r) \leq B' d^3 (\log r)^{1+\alpha}
\]
for every $r\geq d+1$ and $\alpha \in [0,1]$. Since this holds with $\alpha=1$ by \eqref{eq:weak_growth_restate}, it suffices to prove that if this estimate holds with some $\alpha\in (0,1]$ then it holds with $\alpha/2$ also.
Suppose therefore that this estimate holds for some $\alpha\in [0,1]$, and let $r\geq d+1$. If $\log \operatorname{Gr}(Adr^2) \leq \log r$ then we have trivially that $\log \operatorname{Gr}(r) \leq \log \operatorname{Gr}(Adr^2) \leq \log r \leq B' d^3 (\log r)^{1+\alpha/2}$ as required. Otherwise, $\log \operatorname{Gr}(Adr^2) \geq \log r$ and we can absorb the second term into the first in the inequality \eqref{dscsacsac} to read
%
  %
  %
  \[
  \log\Gr(r)\leq 2Cd^{3/2}\sqrt{\log r\log\Gr(Adr^2)}.
  \]
  Substituting in our assumed bound $\log\Gr(r) \leq B' d^3 (\log r)^{1+\alpha}$, we obtain that
\begin{equation*}
\log\Gr(r)
  \leq 2C d^{3/2} \sqrt{ \log r  B'd^3 (\log Ad r^2)^{1+\alpha}}
 \leq B' d^3 (\log r)^{1+\alpha/2}
%
\end{equation*}
by definition of $B'$. This completes the proof that the bound with exponent $\alpha$ implies the bound with exponent $\alpha/2$, and the claim follows by taking $\alpha\downarrow 0$.
\end{proof}

\subsection{Consequences for diffusivity}
\label{subsec:consequences_for_structure_and_diffusivity}

In this section we deduce 
\cref{cor:diffusivity} from \cref{thm:growth_main}.  As mentioned in the introduction, this deduction will invoke the finitary structure theory of transitive graphs of polynomial growth due to Breuillard, Green, and Tao \cite{bgt12} and Tessera and Tointon \cite{TesseraTointon}. More concretely, we will use this theory to deduce that bounded-degree transitive graphs of non-negative Ollivier--Ricci curvature are uniformly \emph{doubling} and then use this doubling property to deduce the diffusive estimate from \cref{prop:displacement_doubling}.
Here we recall that a metric space is said to be \textbf{$k$-doubling} if every bounded subset $S$ of the space can be covered by subsets $S_1,\ldots S_k$  with $\diam(S_i)\leq \diam(S)/2$ for all $1\leq i \leq k$.

\begin{theorem}
\label{thm:structure_doubling}
Let $G$ be a transitive graph satisfying the growth bound $\operatorname{Gr}(r)\leq M r^M$ for all $r\geq 1$ and some $M<\infty$. Then there exists a constant $C=C(M)$ such that $G$ is $C(M)$-doubling.
\end{theorem}

This theorem has the following corollary together with \cref{thm:growth_main}.

\begin{corollary}
For each $d\geq 1$ there exists a constant $C=C(d)$ such that if $G=(V,E)$ is a transitive graph of non-negative Ollivier--Ricci curvature with degrees bounded by $d$ then $G$ is $C(d)$-doubling.
\end{corollary}

\begin{proof}[Proof of \cref{thm:structure_doubling}]
It follows from the finitary structure theory of transitive graphs of polynomial growth \cite[Corollary 1.4]{TesseraTointon} that for each $K<\infty$ there exists $R_0(K)$ and $C_1(K)$ such that if a transitive graph $G$ has 
$\operatorname{Gr}(3r) \leq K\operatorname{Gr}(r)$
for some $r\geq R_0(K)$ then
\[
  \frac{\operatorname{Gr}(R_2)}{\operatorname{Gr}(R_1)} \leq C_1(K) \left(\frac{R_2}{R_1}\right)^{C_1(K)}
\]
for every $R_2\geq R_1 \geq r$.
Letting $K= 3^{2M}$, we observe that if $G$ is a transitive graph satisfying $\operatorname{Gr}(r)\leq M r^M$ for all $r\geq 1$ then
\[
 \frac{1}{k} \sum_{i=1}^k \log_3 \frac{\operatorname{Gr}(3^i R_0(K))}{\operatorname{Gr}(3^{i-1} R_0(K))} \leq \frac{1}{k}\log_3 \operatorname{Gr}(3^k R_0(K)) \leq M + \frac{1}{k} M \log_3 R_0(K) + \frac{1}{k} \log_3 M \leq 2 M
\]
for every $k\geq k_0(M) = M \log_3 R_0(K)+ \log_3M$, so that there must exist $R_0(K) \leq r \leq 3^{k_0(M)}R_0(K)$ with $\operatorname{Gr}(3r)/\operatorname{Gr}(r) \leq K$. By the result recalled above, there exist constants $C_2=C_2(M)$ such that
\begin{equation}
\label{eq:doubling_ratio_growth}
  \frac{\operatorname{Gr}(R_2)}{\operatorname{Gr}(R_1)} \leq C_2 \left(\frac{R_2}{R_1}\right)^{C_2}
\end{equation}
for every $R_2\geq R_1 \geq 3^{k_0(M)}R_0(K)$. Since $\operatorname{Gr}(R_1)$ is bounded by a constant depending only on $M$ in the complementary case $R_1 \leq 3^{k_0(M)}R_0(K)$, we can increase the constant $C_2$ (in a way depending only on $M$) such that the same estimate holds for all $R_2\geq R_1 \geq 1$.

We now prove that this implies bounded doubling constant. Note that if $G=(V,E)$ is a locally finite transitive graph and $S$ is a set of diameter $r$ then, by the Vitali covering lemma, there exist points $x_1,\ldots,x_k$ such that $\bigcup_i B(x_i,r/4) \supseteq S$ and the balls $B(x_i,r/12)$ are disjoint. Since $S$ is contained in a ball of radius $r$, it follows that
$k \operatorname{Gr}(r/12) \leq \operatorname{Gr}(5r/4)$
and hence that 
$k \leq \operatorname{Gr}(5r/4)/\operatorname{Gr}(r/12)$. This and \eqref{eq:doubling_ratio_growth} yield the claimed bound on the doubling constant since $r\geq 1$ was arbitrary and the balls $B(x_i,r/4)$ have diameter at most $r/2$.
\end{proof}

\begin{proof}[Proof of \cref{cor:diffusivity}]
It follows from \eqref{eq:doubling_ratio_growth} and \cref{prop:displacement_doubling} that there exists a constant $C_1(d)$ such that $D_n(x)^2 \leq C_1(d) n$ for every $x\in V$ and $n\geq 1$, and the claim follows since \cref{prop:Lipschitz_concentration} applied with $f=d(x,\,\cdot\,)$ implies that $\mathbf{E}_x d(X_0,X_n)^2 \preceq D_n(x)^2+n$. 
 (Alternatively, it follows from the results of Ding, Lee, and Peres \cite[Theorem 1.5]{ding2013markov} on Markov type-$2$ for doubling metric spaces that if $G$ is any transitive graph that is $\lambda$-doubling then the random walk on $G$ satisfies a diffusive estimate of the form $\mathbf{E}_x d(X_0,X_n)^2 \leq C (\lambda) n$ with $C(\lambda)=O((\log \lambda)^2)$.)
\end{proof}

\subsection{Finitary version}\label{finitaryproof}

As we have already remarked, our results with \cite{CamilloElia} imply that a finitely generated group having a Cayley graph with non-negative Ollivier--Ricci curvature is virtually abelian. We now prove the finitary version of this statement as stated in \Cref{finitary}.


\begin{proof}[Proof of \Cref{finitary}]
Let $S=\{s_1,\dots,s_d\}$ for $d\in\N$. Consider $F_d$, the free group generated by $d$ elements (which we call $s_1,\dots,s_d$) and the natural map $\Psi:F_d\rightarrow G$. Set $N\defeq\ker(\Psi)\unlhd F_d$. Let 
\begin{equation}
    N_7\defeq \{g\in N: |g|\le 7\}\,,
\end{equation}
where $|g|$ denotes the length of the word in the free group. Let also $\tilde N_7$ be the smallest normal subgroup of $F_d$ containing $N_7$. Notice that
$N_7\subseteq \tilde N_7\unlhd N$. We define $\tilde G\defeq {F_d}/{\tilde N_7}$, which is  generated by $S$ (with a slight abuse). 
Notice that 
\begin{equation}
    \{g\in \tilde N_7: |g|\le 7\}=  \{g\in  N: |g|\le 7\}\,.
\end{equation}
Hence, $\operatorname{Cay}(\tilde G,S)$ has  balls of radius $2$ isomorphic to  balls of radius $2$  in $\operatorname{Cay}( G,S)$. In particular, $\operatorname{Cay}(\tilde G,S)$ has  non-negative Ollivier--Ricci curvature so that
%
%
 $\tilde G$ is virtually abelian by \cref{cor:virtually_abelian}. In particular, there exists $\tilde G_N\unlhd \tilde G$ abelian of finite index. 
As $\tilde N_7\unlhd N$, we have a group epimorphism  $\Tilde{G}\rightarrow G$ induced by the projection $F_d\rightarrow G$. In particular, there exists an abelian subgroup of $G$ of index bounded by $|{\tilde G}/{\tilde G_N}|$.
%
%
The claim follows since there are only finitely many possibilities for $N_7$ for $d$ fixed.
\end{proof}

\section{Proof of the main theorem for finite graphs}

We now prove our main theorem for (not-necessarily-transitive) finite graphs, \cref{thm:main_finite}. Much of the proof is the same as that of \cref{prop:log_squared_transitive} and we focus on what needs to be changed.
We now fix for the remainder of the section a finite graph $G=(V,E)$ of non-negative Ollivier--Ricci curvature with maximum degree $d$ and maximum degree ratio $C_\mathrm{rev}= \max_{x,y\in V}\frac{\deg(x)}{\deg(y)}$.

For each vertex $x\in V$ and time $n\geq 1$, we write $D_n(x)=\mathbf{E}_x d(x,X_n)$ for the expected displacement of the $n$-step random walk started at $x$ and write $H_n(x)=-\sum_{y\in V} P^n(x,y)\log P^n(x,y)$ for the Shannon entropy of the $n$-step lazy random walk started at the vertex $x$. We also write $D_n=\frac{1}{|V|}\sum_{x\in V}D_n(x)$ and $H_n=\frac{1}{|V|}\sum_{x\in V} H_n(x)$ for the averages of these quantities over $x\in V$.
%
%
Several of the estimates used in the proof of \cref{prop:log_squared_transitive} adapt unproblematically to (non-transitive) finite graphs, except that many constants now depend on the maximum \emph{degree ratio} $C_\mathrm{rev}$ as well as the maximum degree:
\begin{enumerate}
  \item The probability kernels $P^n$ are $C_\mathrm{rev}$-reversible. Thus, it follows from \cref{prop:entropy_and_cell_size} that if we take $[x]_n$ to be the cell of $x$ in the grand Poisson cell decomposition associated to the probability kernel $P^n$ and define $a_n=\frac{1}{|V|}\sum_{x\in V}\E\log \#[x]_n$ then
  \begin{equation}
  \label{eq:entropy_cell_size_finite_restate}
 \frac{1}{8C_\mathrm{rev}^3}H_n - \frac{1}{2}\leq a_n \leq C_\mathrm{rev}H_n+2\log C_\mathrm{rev}.
\end{equation}
  \item We will apply  \cref{thm:munch} in the form
\begin{equation}
\label{eq:Munch_rearranged_finite_restate}
\log |\Lambda| \leq 4C_\mathrm{rev} \diam(\Lambda) \frac{|\partial W|}{|W|} + \log |W|+\log 2C_\mathrm{rev}.
\end{equation}
  \item The bounds
  \[
 \mathbf{E}_x d(x,X_n)^2 \leq 2n H_n(x) + n \log 4C_\mathrm{rev} \quad \text{ and hence } \quad D_n(x) \leq \sqrt{2n H_n(x)+n \log 4C_\mathrm{rev}}
\]
hold on any bounded-degree graph by the same application of Varopoulos--Carne as the proof of \cref{lem:D_n_upper}. Averaging over $x\in V$, it follows by Jensen that the inequality
 \begin{equation}
 \label{eq:displacement_vs_entropy_finite}
  D_n \leq \sqrt{2n H_n+n \log 4C_\mathrm{rev}} 
\end{equation}
holds on any finite graph for every $n\geq 1$.
\item By the same argument as in the proof of \cref{lem:Gr_bounded_above_by_hat_a}, there exists a universal constant $A$ such that
\[
  \log \#B(x,r) \leq 2 \E\log\#[x]_n + \log 2
\]
for every $x\in V$ and $n,r\geq 1$ with $n\geq A d r^2$. In particular,
\begin{equation}
\label{eq:growth_a_n_finite}
\frac{1}{|V|}\sum_{x\in V} \log \# B(x,r) \leq 2 a_n + \log 2
\end{equation}
whenever $n\geq A d r^2$.
\end{enumerate}

As in the proof of \cref{prop:log_squared_transitive} we will want to work with a truncated version of the grand Poisson cell decomposition associated to a version of $P^n$. However, we can no longer use the same truncation scale everywhere, and for reasons related to this it will be convenient to truncate the cells directly rather than to truncate the kernel defining them. 
For each $n\geq 1$ let $\Pi_n$ be the grand Poisson cell decomposition associated to the kernel $P^n$, where we take these cell decompositions to be independent for different values of $n$. Given $\Pi_n$, we say that a vertex $x$ is \textbf{$n$-good} if $d(x,X_n(x)) \leq \min\{D_n(x), D_n(X_n(x))\} + 4 \sqrt{n\log n}$, where $(X_n(x))_{x\in V}$ is the coupled family of lazy $n$-step random walks used to define the partition $\Pi_n$, and say that $x$ is \textbf{$n$-bad} otherwise.

\begin{lemma}
\label{lem:bad}
$\frac{1}{|V|}\sum_{x\in V}\P(x \text{ is $n$-bad}) \preceq C_\mathrm{rev} n^{-2}$ for every $n\geq 1$.
\end{lemma}

\begin{proof}[Proof of \cref{lem:bad}]
The probability that $x$ is $n$-bad is bounded by $\bP_x(d(x,X_n) \geq D_n(x)+4\sqrt{n \log n})+\bP_x(d(x,X_n) \geq D_n(X_n)+4\sqrt{n \log n})$.
Applying \cref{prop:Lipschitz_concentration} with $f=d(x,\,\cdot\,)$ shows that the first of these two probabilities is at most $2n^{-2}$ for every $x\in V$. We can then use the mass-transport principle \eqref{eq:finiteMTP} to compute that
\begin{align*}
  \frac{1}{|V|}\sum_{x\in V}\P(x \text{ is $n$-bad}) 
&\leq 2n^{-2}+\frac{1}{|V|}\sum_{x\in V}\E\left[ \sum_{y\in V} P^n(x,y) \mathbbm{1}(d(x,y) \geq D_n(y)+4\sqrt{n \log n}) \right]
\nonumber\\&
= 2n^{-2}+\frac{1}{|V|}\sum_{x\in V}\E\left[ \sum_{y\in V} P^n(y,x) \mathbbm{1}(d(x,y) \geq D_n(x)+4\sqrt{n \log n}) \right]
\nonumber\\&
= 2n^{-2}+\frac{1}{|V|}\sum_{x\in V}\E\left[ \sum_{y\in V} \frac{\deg(x)}{\deg(y)}P^n(x,y) \mathbbm{1}(d(x,y) \geq D_n(x)+4\sqrt{n \log n}) \right]
  \nonumber\\&\leq 2n^{-2}+C_\mathrm{rev}
   \frac{1}{|V|}\sum_{x\in V}\E\left[\bP_x(d(x,X_n) \geq D_n(x)+4\sqrt{n \log n})\right]
   \preceq C_\mathrm{rev} n^{-2},
\end{align*}
where we used \cref{prop:Lipschitz_concentration} again in the last line and the time-reversal identity $\deg(x)P^n(x,y)=\deg(y)P^n(y,x)$ in the third line. 
\end{proof}

 We define the partition $\widetilde \Pi_n$ by removing every $n$-bad vertex from its cell and placing it in a singleton cell, so that the cell $\widetilde{[x]}_n$ of $x$ in $\widetilde{\Pi}_n$ is equal to the set of $n$-good elements of $[x]_n$ if $x$ is $n$-good and equal to $\{x\}$ otherwise.
Observe that for each $n$-good vertex $x\in V$, the diameter of the cell $\widetilde{[x]}_n$ is bounded by $2D_n(X_n(x))+8\sqrt{n\log n}$, with the location of the random walk $X_n(y)$ being constant over the cell $[x]_n$ by definition of the grand Poisson cell decomposition, and \cref{lem:D_n_Lipschitz} implies that this diameter  is bounded by $6D_n(x)+24\sqrt{n\log n}$ when $x$ is $n$-good.

\medskip

The next lemma shows that passing from $\Pi_n$ to $\widetilde \Pi_n$ does not significantly change expected log cell sizes.

\begin{lemma} 
\label{lem:bad_log_vol}
It holds that
\[
 \frac{1}{|V|}\sum_{x\in V} \E \log \# [x]_n -\frac{1}{|V|}\sum_{x\in V}\E\log \#\widetilde{[x]}_n\preceq     C_\mathrm{rev} n^{-1}\log d+1
\]
for every $n\geq 1$.
\end{lemma}

\begin{proof}[Proof of \cref{lem:bad_log_vol}]
We can compute that
\begin{align}
 \frac{1}{|V|}\sum_{x\in V}\P\left(\text{$x$ is $n$-good and $\#\widetilde{[x]}_n\leq \frac{1}{2}\#[x]_n$}\right) 
&\leq \frac{2}{|V|}\sum_{x\in V}\E \frac{1}{\#[x]_n} \sum_{y\in [x]_n} \mathbbm{1}(\text{$x$ is $n$-good, $y$ is $n$-bad})
 \nonumber\\& = \frac{2}{|V|}\sum_{x\in V} \E \frac{1}{\#[x]_n} \sum_{y\in [x]_n} \mathbbm{1}(\text{$x$ is $n$-bad, $y$ is $n$-good})
 \nonumber\\ &\leq \frac{2}{|V|}\sum_{x\in V} \P(\text{$x$ is $n$-bad}) \preceq C_\mathrm{rev} n^{-2},
 \label{eq:bad_cell_comparison}
\end{align}
where the first line follows by Markov's inequality, the second follows from the mass-transport principle, and the final inequality follows from \cref{lem:bad}. To apply this inequality, we observe that
\begin{multline*}
   \frac{1}{|V|}\sum_{x\in V} \E\log \#\widetilde{[x]}_n \geq  \frac{1}{|V|}\sum_{x\in V}\E \left[\mathbbm{1}(\text{$x$ is $n$-good and $\#\widetilde{[x]}_n\geq \frac{1}{2}\#[x]_n$})\log \frac{1}{2}\#[x]_n \right]
    \\\geq \frac{1}{|V|}\sum_{x\in V}\left(\E \log \frac{1}{2}\#[x]_n - \E \left[\mathbbm{1}(\text{$x$ is $n$-bad or $x$ is $n$-good and $\#\widetilde{[x]}_n\leq \frac{1}{2}\#[x]_n$})\log \frac{1}{2}\#[x]_n \right]\right).
\end{multline*}
Using the trivial bound $\log \#[x]_n \preceq (n+1)\log (d+1)$, which holds because the cell $[x]_n$ is contained within  the ball of radius $n$ around $X_n(x)$, we can substitute \eqref{eq:bad_cell_comparison} into this inequality and rearrange to yield the claim.
\end{proof}
\medskip

We now derive our main functional inequality for the sequence $a_n\defeq \frac{1}{|V|}\sum_{x\in V}\E \log \#[x]_n$.

\begin{prop}
\label{lem:a_n_functional_inequality_with_Dn_finite} The inequality
\[
  a_n-a_m \preceq \frac{ d^{3/2}C_\mathrm{rev}}{\sqrt{m}}(D_n+\sqrt{n\log n})+ C_\mathrm{rev} + \frac{d C_\mathrm{rev}^2}{n}
\]
holds for every $n>m\geq 1$.
\end{prop}


\begin{proof}[Proof of \cref{lem:a_n_functional_inequality_with_Dn_finite}]
For each $m<n$ let
 $\widetilde{[x]}_{n,m}= \widetilde{[x]}_{n} \cap [x]_{m}$ where we take the two different grand Poisson cell decompositions associated to $ P^n$ and $ P^m$ to be independent. (Note that we do not throw away $m$-bad points in the partition $\widetilde \Pi_m$.) It follows from \eqref{eq:Munch_rearranged_finite_restate} that
\[
  \log \#\widetilde{[x]}_n \leq \log \#\widetilde{[x]}_{n,m} + \log 2 C_\mathrm{rev} + 4 C_\mathrm{rev} \diam (\widetilde{[x]}_n) \frac{\#\partial \widetilde{[x]}_{n,m}}{\# \widetilde{[x]}_{n,m}}.
\]
Bounding $\log \#\widetilde{[x]}_{n,m} \leq \log \#[x]_{m}$, taking expectations and averaging over $x$ yields that
\[
 \frac{1}{|V|}\sum_{x\in V}\E \log \#\widetilde{[x]}_n  -\frac{1}{|V|}\sum_{x\in V} \E  \log \#[x]_m \preceq C_\mathrm{rev} 
  +  \frac{C_\mathrm{rev}}{|V|}\sum_{x\in V} \E\left[\diam (\widetilde{[x]}_n) \frac{\#\partial \widetilde{[x]}_{n,m}}{\# \widetilde{[x]}_{n,m}}\right]
\]
and hence by \cref{lem:bad_log_vol} that
\begin{equation}
\label{eq:finite_Munch_expectation}
a_n-a_m \preceq C_\mathrm{rev} 
  \left(1+\frac{\log d}{n}\right) 
  +  \frac{C_\mathrm{rev}}{|V|}\sum_{x\in V} \E\left[\diam (\widetilde{[x]}_n) \frac{\#\partial \widetilde{[x]}_{n,m}}{\# \widetilde{[x]}_{n,m}}\right].
\end{equation}
Since $\diam (\widetilde{[x]}_n)$ is constant on $\widetilde{[x]}_{n,m}$, we can apply \cref{lem:cell_isoperimetry_finite} to obtain that
\begin{equation}
  \frac{1}{|V|} \sum_{x\in V}\E\left[\diam (\widetilde{[x]}_n) \frac{\#\partial \widetilde{[x]}_{n,m}}{\# \widetilde{[x]}_{n,m}}\right] = \frac{1}{|V|} \sum_{x\in V}  \E\left[\diam (\widetilde{[x]}_n)\sum_{z\sim x}\mathbbm{1}(z\notin \widetilde{[x]}_{n,m})\right].
\label{eq:finite_diam_boundary_rewrite}
\end{equation}
Observe that if $z\notin \widetilde{[x]}_{n,m}$ then either at least one of $x$ or $z$ is $n$-bad, or $z\notin [x]_n$, or $z\notin [x]_m$.
The mass-transport principle and \cref{lem:bad} yield that
\[
  \frac{1}{|V|} \sum_{x\in V} \sum_{z\sim x} \P(z \text{ is $n$-bad}) = \frac{1}{|V|} \sum_{x\in V} \sum_{z\sim x} \P(x \text{ is $n$-bad}) \preceq \frac{d C_\mathrm{rev}}{n^2},
\]
and using the trivial bound $\diam (\widetilde{[x]}_n)\leq 2n$ yields that
\begin{equation}
\frac{1}{|V|} \sum_{x\in V}  \E\left[\diam (\widetilde{[x]}_n)\sum_{z\sim x}\mathbbm{1}(\text{$x$ or $z$ is $n$-bad})\right] \preceq \frac{d C_\mathrm{rev}}{n}.
\label{eq:bad_diameter_boundary}
\end{equation}
On the other hand, we recall the diameter bound $\diam (\widetilde{[x]}_n) \leq 6 D_n(x) + 24 \sqrt{n \log n }$ for $x$ $n$-good.
Since $D_n(x)$ is not random, we can take it out of the expectation to obtain that
\begin{multline}
\frac{1}{|V|} \sum_{x\in V}  \E\left[\diam (\widetilde{[x]}_n)\sum_{z\sim x}\mathbbm{1}(x,z \text{ $n$-good, } z\notin [x]_{n} \cap [x]_m)\right]
\\
\leq \frac{1}{|V|} \sum_{x\in V} (6D_n(x)+24 \sqrt{n\log n}) \sum_{z\sim x} \left[\P(z\notin [x]_n) + \P(z\notin [x]_m)\right] 
 \preceq \frac{d^{3/2}}{\sqrt{m}} (D_n+ \sqrt{n \log n}),
 \label{eq:good_diameter_boundary}
\end{multline}
where we used \eqref{eq:coupling_TV} and \Cref{thm:Salez} in the last line.
The claim follows from \eqref{eq:finite_Munch_expectation}, \eqref{eq:finite_diam_boundary_rewrite}, \eqref{eq:bad_diameter_boundary}, and \eqref{eq:good_diameter_boundary}. \qedhere
\end{proof}

\begin{prop}\label{prop:vcdscsd}
    It holds that 
    \begin{equation}\notag
        \frac{1}{|V|}\sum_{x\in V} \log \#B(x,r)\preceq  d^3 C_\mathrm{rev}^5 (\log r)^2
    \end{equation}
    for every $r\ge d+1$.
\end{prop}
\begin{proof}[Proof of \cref{prop:vcdscsd}]
The lower bound of \eqref{eq:entropy_cell_size_finite_restate} can be rearranged to yield that
$H_n \leq 4 C_\mathrm{rev}^3 (2a_n+1)$, and substituting this bound into \eqref{eq:displacement_vs_entropy_finite} yields that
\[
  D_n \leq \sqrt{8n C_\mathrm{rev}^3 (2a_n+1)+n \log 4 C_\mathrm{rev}} \preceq \sqrt{C_\mathrm{rev}^3 n(a_n+1)}.
\]
Substituting this bound into \cref{lem:a_n_functional_inequality_with_Dn_finite} yields that
\[
  a_n -a_m \preceq \frac{ d^{3/2}C_\mathrm{rev}}{m^{1/2}}(\sqrt{C_\mathrm{rev}^3 n(a_n+1)}+\sqrt{n\log n})+ C_\mathrm{rev} + \frac{d C_\mathrm{rev}^2}{n}
\]
and hence that
\begin{equation}
\label{eq:displayed_inequality}
  a_n -a_m \preceq \sqrt{\frac{d^3 C_\mathrm{rev}^5 n a_n}{m}} 
\end{equation}
whenever $n> m \geq d+1$ and $a_n \geq \log n$.
It now follows by the same argument as the proof of \cref{prop:log_squared_transitive} that $a_{d2^k}\preceq d^3 C_\mathrm{rev}^5 k^2$ for every $k\geq 1$, and hence, applying \eqref{eq:displayed_inequality} once more with $m=d2^k\leq n\leq d2^{k+1}$ (the case $a_n\leq \log n$ being trivial), that $a_n \preceq d^3C_\mathrm{rev}^5(\log n)^2$ for every $n\geq 2$.  The conclusion then follows from \eqref{eq:growth_a_n_finite}.
\end{proof}

\begin{proof}[Proof of \cref{thm:main_finite}]
Fix for the moment $x\in V$. As in the proof of \Cref{prop:main_transitive_functional_inequality},  
we have
\begin{equation}\notag
   \sum_{\ell=1}^{k}\frac{D_{2^\ell}(x)}{2^{\ell/2}} \preceq k+k(\log C_{\mathrm{rev}})^{1/2}+\sqrt{k\log \#B(x,2^{k+4})}
\end{equation}
for every $k\ge 1$. Averaging over $x\in V$ and using Cauchy--Schwarz, we have
\begin{equation}\notag
    \sum_{\ell=1}^{k}\frac{D_{2^\ell}}{2^{\ell/2}} \preceq k+k(\log C_{\mathrm{rev}})^{1/2}+\sqrt{k\frac{1}{|V|}\sum_{x\in V} \log \#B(x,2^{k+4})}
\end{equation}
for every $k\ge 1$.
Using (as in the proof of \Cref{prop:main_transitive_functional_inequality}) \Cref{lem:a_n_functional_inequality_with_Dn_finite}, we obtain
\begin{align*}
    a_{2^k}&\preceq a_2+\sum_{\ell=1}^k \bigg(d^{3/2}C_{\mathrm{rev}}\frac{D_{2^\ell}}{2^{\ell/2}}+d^{3/2}C_{\mathrm{rev}}\ell^{1/2}+ C_\mathrm{rev} + \frac{d C_\mathrm{rev}^2}{2^\ell}\bigg)\\
    &\preceq d^{3/2}C_{\mathrm{rev}}k^{3/2}+kC_{\mathrm{rev}}+dC_{\mathrm{rev}}^2+d^{3/2}C_{\mathrm{rev}}\sum_{\ell=1}^k \frac{D_{2^\ell}}{2^{\ell/2}}\\
    &\preceq d^{3/2}C_{\mathrm{rev}}k^{3/2}+dC_{\mathrm{rev}}^2+d^{3/2}C_{\mathrm{rev}}k(\log C_{\mathrm{rev}})^{1/2}\\
    &\qquad\qquad+d^{3/2}C_{\mathrm{rev}}\sqrt{k\frac{1}{|V|}\sum_{x\in V} \log \#B(x,2^{k+4})}.
\end{align*}
Following again the proof of \Cref{prop:main_transitive_functional_inequality} and using \eqref{eq:growth_a_n_finite} and \Cref{prop:vcdscsd}, we obtain
\begin{align*}
  \frac{1}{|V|}\sum_{x\in V} \log \# B(x,r)&\preceq  d^{3/2}C_{\mathrm{rev}}(\log r)^{3/2}+dC_{\mathrm{rev}}^2+d^{3/2}C_{\mathrm{rev}}\log r(\log C_{\mathrm{rev}})^{1/2}\\
 &\qquad\qquad+d^{3/2}C_{\mathrm{rev}}(\log r)^{1/2}d^{3/2}C_{\mathrm{rev}}^{5/2}\log r
\end{align*}
for every $r\ge d+1$. The case $r \leq d$  holds trivially.

Finally, the claimed bound on the displacement follows since
\[
  \frac{1}{|V|}\sum_{x\in V}\mathbf{E}_x d(x,X_n)^2 \preceq nH_n + n \log 4 C_\mathrm{rev} \preceq  n C_\mathrm{rev}^3 (a_n+1) \preceq n C_\mathrm{rev}^3 \frac{1}{|V|}\sum_{x\in V}\log \#B(x,2n),
\]
for every $n\geq 2$.
\end{proof}

\subsection*{Acknowledgements}
TH was supported by NSF grant DMS-1928930 and a Packard Fellowship for Science and Engineering. Part of this work was carried out while CB was a postdoc at the Institute for Advanced Study, Princeton, supported by the NSF  grant  DMS-2424441.  We thank Elia Bruè, Kunal Chawla, Antoine Song, and Tianyi Zheng for helpful comments and discussions, and give special thanks to Jonathan Hermon for many stimulating discussions that catalyzed the project and for carefully reading an earlier draft of the paper.

\footnotesize{
\bibliographystyle{abbrv}
\bibliography{big_bib_file}
}

\normalsize
\begin{appendix}
\section{Proof of the isoperimetric inequality}
\label{sec:Munch_appendix}

In this appendix we prove \cref{thm:munch}. We stress that the primary purpose of this section is expository, and that the proof we give here does not significantly differ in mathematical content from that of \cite[Theorem 6.1]{munch2023ollivier}. We also note that some key parts of the proof, including the Laplacian separation principle \cref{thm:Laplacian_Separation}, originate in the earlier work of Hua and the third named author \cite{hua2025every}. See also \cite[Section 4.3]{li2024convergence} for an alternative approach to the Laplacian separation principle.

\subsection{The Laplacian separation principle}

Given a connected, locally finite graph $G=(V,E)$, we write $\Lip(1)$ for the space of $1$-Lipschitz functions on $G$, noting that $\{f\in \Lip(1):f(x_0)=0\}$ is compact in the product topology for each $x_0\in V$.
Throughout this section we will let $P$ denote the Markov operator for the lazy random walk (indeed, all walks considered will be lazy) and use the sign convention 
$\Delta=P-I$ for the Laplacian.
Note that if $\theta \in [0,1]$ then
\[
  f+\theta \Delta f = \theta P f+(1-\theta)f
\]
is a convex combination of $f$ and $Pf$. 
As such, if $G$ has non-negative Ollivier--Ricci curvature and $f$ is $1$-Lipschitz then $f+\theta \Delta f$ is a convex combination of two $1$-Lipschitz functions and is therefore $1$-Lipschitz also.

\begin{theorem}[Laplacian Separation Principle]
\label{thm:Laplacian_Separation}
Let $G=(V,E)$ be a connected, locally finite graph with non-negative Ollivier--Ricci curvature. Suppose $V$ is written as a disjoint union $V = X \cup K \cup Y$ for a finite non-empty connected set $K$ such that there are no edges between $X$ and $Y$. There exists a function $f \in \operatorname{Lip}(1)$ and a constant $\lambda\in \R$ such that
\begin{enumerate}
    \item $\Delta f = \lambda$ on $K$.
    \item $f(x) = \min\{g(x):g \in \operatorname{Lip}(1): g|_{K} = f|_{K}\}$ and $\Delta f(x) \geq \lambda$ for every $x\in X$.
    \item $f(y) = \max\{g(y):g \in \operatorname{Lip}(1): g|_{K} = f|_{K}\}$ and $\Delta f(y) \leq \lambda$ for every $y\in Y$.
\end{enumerate}
\end{theorem}

Note that the sets $X$ and $Y$ are \emph{not} required to be non-empty, and indeed we will later apply the theorem with $X=\emptyset$ via the following corollary.

\begin{corollary}
\label{cor:Laplacian_Separation}
Let $G=(V,E)$ be a connected, locally finite graph with non-negative Ollivier--Ricci curvature and suppose that $K\subseteq V$ is a non-empty finite connected set. There exists a function $f \in \operatorname{Lip}(1)$ and a constant $\lambda\in \R$ such that
\begin{enumerate}
  \item $\Delta f = \lambda$ on $K$ and $\Delta f \leq \lambda$ on $V\setminus K$.
  \item For every $x\in V\setminus K$ there exists $y\sim x$ with $f(y) \leq f(x)-1$.
  \item $f$ attains its global minimum on $K$.
\end{enumerate}
\end{corollary}

\medskip

Before proving \cref{thm:Laplacian_Separation,cor:Laplacian_Separation}, let us first introduce the Lipschitz extension map and its basic properties. For the remainder of the section we fix the graph $G=(V,E)$ and the decomposition $V=X\cup K \cup Y$ as in \cref{thm:Laplacian_Separation}.
Recall that the set $\operatorname{Lip}(1)$ is closed under convex combinations, maxima, and minima.
Let $V=X\cup K \cup Y$ be written as a disjoint union above and for each function $f\in \operatorname{Lip}(1)$ define
\[
  Sf(x) = \begin{cases}
  f(x)&x\in K \\
  \min\{g(x):g\in \operatorname{Lip}(1),\, g|_K= f|_K\} &x\in X \\
  \max\{g(x):g\in \operatorname{Lip}(1),\, g|_K= f|_K\} &x\in Y
  \end{cases}.
\]
$S$ defines an idempotent function mapping $\operatorname{Lip}(1)$ to $\operatorname{Lip}(1)$. Observe that if $f\in \operatorname{Lip}(1)$ then we can express $Sf(x)$ equivalently as
\[
  Sf(x) = \begin{cases}
  f(x)&x\in K \\
  \min\{g(x):g\in \operatorname{Lip}(1),\, g|_K\geq f|_K\} &x\in X \\
  \max\{g(x):g\in \operatorname{Lip}(1),\, g|_K\leq f|_K\} &x\in Y
  \end{cases},
\]
where the optimization problems are now stated with inequalities rather than equalities, so that if $f_1,f_2 \in \operatorname{Lip}(1)$ are such that $f_1 \geq f_2$ on $K$ then $S(f_1) \geq S(f_2)$ on the entire set $V$; we refer to this as the \textbf{order-preserving} property of $S$. Note also that if $c$ is a constant then $S(f+c)=S(f)+c$. Note that we can also write $Sf$ explicitly as
\begin{equation}
\label{eq:S_explicit}
  Sf(x)=\begin{cases}f(x)&x\in K\\
  \max\{f(z)-d(x,z):z\in K\} & x \in X\\
  \min\{f(z)+d(x,z):z\in K\} & x \in Y
  \end{cases}.
\end{equation}
(This representation makes it obvious that the map $f\mapsto S(f)$ is continuous on $\operatorname{Lip}(1)$ equipped with the product topology.)
We define $\cF_0=S(\Lip(1))=\{f\in \Lip(1) : Sf=f\}$.

\begin{lemma}
\label{lem:Laplacian_on_X_and_Y}
If $f\in \cF_0$ then $\Delta f(x)\geq \min_K \Delta f$ for all $x\in X$ and $\Delta f(x) \leq \max_K \Delta f$ for all $x\in Y$.
\end{lemma}

\begin{proof}[Proof of \cref{lem:Laplacian_on_X_and_Y}]
We prove that if $f\in \cF_0$ then $\Delta f(x)\geq \min_K \Delta f$ for all $x\in X$, the corresponding estimate for $Y$ being similar. Let $C= \min_K \Delta f$.
By non-negative curvature the function  $f+\eps \Delta f-\eps C$ is $1$-Lipschitz  for each $0\leq \eps\leq 1$, and satisfies $f+\eps \Delta f-\eps C \geq f$ on $K$ by choice of $C$.
 Since $Sf=f$, meaning that $f|_X$ is smaller than every other $1$-Lipschitz function whose restriction to $K$ is larger than $f$, we must have that
$f(x) \leq f(x)+\eps \Delta f(x)-\eps C$ for every $x\in X$.
This inequality rearranges to give the desired inequality for $\Delta f$.
\end{proof}

Let $\cF_1$ be the set of functions in $\cF_0$ maximizing $\min_K \Delta f$ (this minimum being attained since 
$K$ is finite and since the maximal value of $\min_K \Delta f$ is unchanged if we restrict to the compact set $\{f\in \cF_0:f(x_0)=0\}$ for some base point $x_0\in V$). Our next lemma establishes a stability property of the set $\cF_1$ under the map $f\mapsto S(f+\eps \Delta f)$ for small $\eps>0$ and shows that this map reduces the size of the set of points in $K$ where the minimum $\min_K \Delta f$ is attained.

\begin{lemma}
\label{lem:cF_2_property}
For each $f\in \cF_1$
there exists $\eps_0=\eps_0(f)>0$ such that if $0\leq \eps \leq \eps_0$ then
 $S(f+\eps\Delta f)\in \cF_1$ and
 \begin{equation}\label{brtvfdv}
     \{x\in K : \Delta S(f+\eps\Delta f)(x) = \min_K \Delta S(f+\eps\Delta f)\}
\subseteq 
\{x\in K : \Delta f(x) = \min_K \Delta f\}.
 \end{equation}
\end{lemma}

\begin{proof}[Proof of \cref{lem:cF_2_property}]
 Let $f\in \cF_1$, let $\lambda = \min_{K}\Delta f$, and let $K'=\{x\in K:\Delta f(x)=\lambda\}$.
Since $S$ is continuous and $S(f)=f$, there exists $\eps_0\in (0,1)$ such that
   $g=g_\eps=S(f+\eps \Delta f)$ satisfies $\Delta g(x) > \lambda$ for every $0\leq \eps\leq \eps_0$ and $x\in K \setminus K'$ (this claim holding vacuously if $K\setminus K'$ is empty). 
Since $g\in \cF_0$ and $f\in \cF_1$, we must have by definition of $\cF_1$ that
\begin{equation}
\label{eq:minDeltag1}
\min_K \Delta g \leq \min_K \Delta f = \lambda
\end{equation}
and hence that
\begin{equation}
\label{eq:Deltag_min_set}
  \{x\in K : \Delta g(x)=\min_K \Delta g\} \subseteq K'.
\end{equation}
   Since $\eps \Delta f \geq \eps \lambda$ on $K$ and $S(f)=f$, it follows by the order-preserving property of $S$ that
   \[
     g =S(f+\eps \Delta f) \geq S(f+\eps \lambda)=S(f)+\eps \lambda=f+\eps \lambda
   \]
   on $V$. This inequality is an equality on $K'$ since on this set we have that $g = f+\eps \Delta f = f+\eps \lambda$. As such, for each $x\in K'$ we must have that 
   \begin{equation}
\label{eq:Deltag_Deltaf_comparison_on_K'}
   \Delta g(x) = \sum_{y}P(x,y)(g(y)-g(x)) \geq \sum_{y}P(x,y)(f(y)+\eps \lambda-(f(x)+\eps \lambda)) = \Delta f(x),
   \end{equation}
   and in light of \eqref{eq:minDeltag1} and \eqref{eq:Deltag_min_set} we must in fact have that $\min_K \Delta g=\lambda$, so that $g\in \cF_1$ as required  and \eqref{brtvfdv} follows as well.
\end{proof}





We now apply \cref{lem:cF_2_property} to prove that $\Delta f$ is constant on $K$ when $f\in \cF_1$. It is in this step that we use the assumption that $K$ is connected.

\begin{lemma}
\label{lem:DeltaF1_flat}
If $f\in \cF_1$ then $\Delta f|_K \equiv \max_{\cF_0} \min_K \Delta f$.
\end{lemma}

\begin{proof}[Proof of \cref{lem:DeltaF1_flat}]
Starting with an arbitrary element of $\cF_1$ and applying \cref{lem:cF_2_property} iteratively (at most $|K|$ times), we deduce that there exists $\eps>0$ and a function $f\in \cF_1$ such that
\[
 K'\defeq \{x\in K : \Delta S(f+\eps \Delta f)(x) = \lambda\} = \{x\in K : \Delta f(x) = \lambda\}.
\]
In order for the inequality \eqref{eq:Deltag_Deltaf_comparison_on_K'} not to be strict for $x\in K'$, we must have that $g(y)=f(y)+\eps \lambda$ for every $y \sim x$. If $y$ is both adjacent to $x$ and belongs to $K$ this means that $\Delta f(y)=\lambda$, and since $K$ was assumed to be connected this implies that $K'=K$. Since the sets $\{x\in K:\Delta f(x)=\lambda\}$ are decreasing along the iteration by \cref{lem:cF_2_property}, and the last of them is equal to $K$, the same is true of the first, which was associated to an arbitrary element of $\cF_1$.
\end{proof}

\begin{proof}[Proof of \cref{thm:Laplacian_Separation}]
This is an immediate consequence of \cref{lem:Laplacian_on_X_and_Y,lem:DeltaF1_flat}.
\end{proof}

\begin{proof}[Proof of \cref{cor:Laplacian_Separation}]
This follows by applying \cref{thm:Laplacian_Separation} with $X=\emptyset$ and $Y=V\setminus K$ and noting that 
since $f$ is equal to its maximal $1$-Lipschitz extension to $V\setminus K$ it must satisfy the second and third claimed properties by \eqref{eq:S_explicit}.
\end{proof}

\subsection{Proof of \cref{thm:munch}}

We now apply the Laplacian separation principle (via \cref{cor:Laplacian_Separation}) to prove \cref{thm:munch}. We will require the following lemma, which lets us transform a function with small Laplacian on $X$ into a different function with larger negative Laplacian via composition with a well-chosen concave function.

\begin{lemma}
\label{lem:Laplacian_transformation_phi}
Let $G=(V,E)$ be a connected graph of non-negative Ollivier--Ricci curvature with maximum degree at most $d$.
Let $X\subseteq V$, let $\lambda\in \R$, and let $f:V\to [0,\infty)$ be a $1$-Lipschitz function such that 
\begin{enumerate}
  \item $\Delta f \leq \lambda$ on $X$.
  \item For every $u\in X$ there exists $v\sim u$ with $f(v)\leq f(u)-1$. 
\end{enumerate}
For each $R\geq 1$ there exists a concave increasing function $\phi:[0,\infty)\to [0,\infty)$ with $\phi'\leq 1$ such that $\phi$ is constant on $[R,\infty)$ and
\[
  \Delta (\phi \circ f) \leq 
\frac{-\lambda}{2\exp(4d\lambda R)-2} \mathbbm{1}(f(x)\leq R)
\]
for every $x\in X$, with the right hand side interpreted as $-\mathbbm{1}(f(x)\leq R)/8dR$ if $\lambda=0$.
\end{lemma}

Note that the quantity $-\lambda/(2\exp(4d\lambda R)-2)$ appearing in this bound is always negative.

Before proving this lemma, let us first see how it can be used to prove \cref{thm:munch} in conjunction with \cref{cor:Laplacian_Separation}.

\begin{proof}[Proof of \cref{thm:munch}]
We may assume that $\pi(\Lambda)>2\pi(W)$, the claimed inequality holding trivially otherwise, so that $\Lambda$ has at least two elements and hence $\diam(\Lambda)\geq 1$.
It suffices to prove the claim for connected sets $W$, the general case following by summing over connected components.
Recall that the Dirichlet inner product between two functions $f$ and $g$ is defined by
\[
  \cE(f,g) = \frac{1}{2}\sum_{x,y}\deg(x)P(x,y) (f(y)-f(x))(g(y)-g(x))
  = \langle f,-\Delta g\rangle_\pi,
\]
with the final equality holding whenever all relevant sums converge absolutely and where we recall that $\langle f, g\rangle_\pi=\sum_{x\in V}\deg(x) f(x)g(x)$.
Observe that if $f$ is any $1$-Lipschitz function then
\[
 |\partial W| \geq |\cE(\mathbbm{1}_W,-f)|= |\langle \mathbbm{1}_W,\Delta f \rangle_\pi|.
\]
As such, if $\lambda$ is the constant obtained by applying \cref{cor:Laplacian_Separation} with $K=W$ then 
\begin{equation}
\label{eq:big_lambda}
   |\partial W| \geq |\lambda| \pi(W).
\end{equation}
This already establishes the claim when $|\lambda|$ is sufficiently large. 
We will now use \cref{lem:Laplacian_transformation_phi} to establish an alternative inequality that instead gets stronger for small values of $|\lambda|$.
Let $f$ be the function whose existence is ensured by \cref{cor:Laplacian_Separation} (applied with $K=W$ and $Y=V\setminus K$), which we normalize to have minimum value $0$ (by adding a constant if necessary).
Let $R=\operatorname{diam}(\Lambda)$, so that $f \leq R$ on $\Lambda$ (since $f$ is $1$-Lipschitz and attains its minimum of zero on $W\subseteq \Lambda$). Applying \cref{lem:Laplacian_transformation_phi} (with $X=V\setminus W$) to this function $f$, we obtain a non-negative $1$-Lipschitz function $g=\phi \circ f$ such that
\[
  \Delta g \leq \frac{-\lambda}{2\exp(4\lambda dR)-2}<0
\]
on $\{x\in V \setminus W :f(x) \leq R\}$, $g$ is constant on $\{x\in V: f(x)\geq R\}$, and $\Delta g \leq 0$ on $\{x\in V: f(x)\geq R\}$. Since for every $x\in V\setminus W$ there exists $y\sim x$ with $f(y)\leq f(x)-1$, and since $f\geq 0$, we have $d(x,W)\leq f(x)$ for every $x\in V$, so that the set $\{x\in V : f(x)\leq R+1\}$ is finite and the function $\Delta g$ is finitely supported. As such, we can safely write the Dirichlet inner product $\cE(\mathbbm{1}_{V\setminus W},g)$ in terms of the inner product $\langle \mathbbm{1}_{V\setminus W},-\Delta g \rangle_\pi$ to obtain that
\begin{multline}
  |\partial W| \geq \cE(\mathbbm{1}_{V\setminus W},g) = \langle \mathbbm{1}_{V\setminus W},-\Delta g \rangle_\pi \\\geq \frac{\lambda}{2\exp(4\lambda dR)-2} \pi(\{x\in V \setminus W :f(x) \leq R\}) \geq \frac{\lambda [\pi(\Lambda)-\pi(W)]}{2\exp(4\lambda dR)-2}.
  \label{eq:small_lambda}
\end{multline}
Putting together \eqref{eq:big_lambda} and \eqref{eq:small_lambda} we obtain that
\[
  \frac{|\partial W|}{\pi(W)} \geq \max\left\{|\lambda|,\frac{\lambda}{2\exp(4\lambda dR)-2} \cdot \frac{\pi(\Lambda)-\pi(W)}{\pi(W)} \right\} \eqdef \Psi(\lambda).
\]
Since $\frac{\lambda}{2\exp(4\lambda dR)-2}$ is a decreasing function of $\lambda$, the minimum $\min_\lambda \Psi(\lambda)$ is attained at the positive value of $\lambda$ satisfying
\[
  \lambda= \frac{\lambda}{2\exp(4\lambda dR)-2} \cdot \frac{\pi(\Lambda)-\pi(W)}{\pi(W)},
\]
which we can compute to be
\[
\lambda=\frac{1}{4dR}\log \left[1+\frac{\pi(\Lambda)-\pi(W)}{2\pi(W)}\right].
\]
As such, we have that
\[
  \min_\lambda \Psi(\lambda) = \frac{1}{4dR}\log \left[1+\frac{\pi(\Lambda)-\pi(W)}{2\pi(W)}\right] \geq \frac{1}{4dR}\log \left[\frac{\pi(\Lambda)}{2\pi(W)}\right]
\]
and hence that
\[
  \frac{|\partial W|}{\pi(W)} \geq \frac{1}{4dR}\log \left[\frac{\pi(\Lambda)}{2\pi(W)}\right] 
\]
as claimed.
\end{proof}

It remains to prove \cref{lem:Laplacian_transformation_phi}. We begin with the following inequality version of the second-order Taylor expansion applying to Laplacians of compositions with concave functions.

\begin{lemma}
\label{lem:chain_rule}
Let $G=(V,E)$ be a graph with degrees bounded by $d$, let $f:V\to \R$ be a function taking values in an interval $I$, and suppose that $\phi:I\to \R$ is a concave function. Then
\[
  \Delta (\phi \circ f)(x) \leq \phi'(f(x))\Delta f(x) + \frac{1}{4d} \phi''(f(x)) \cdot (f(x)- \min\{f(y) : y\sim x \text{ or } y=x\} )^2
\]
for every $x\in V$ such that $\phi'''(t)\geq 0$ for every $t\in [\min\{f(y) : y\sim x \text{ or } y=x\},f(x)]$.
\end{lemma}

\begin{proof}[Proof of \cref{lem:chain_rule}]
Let $x\in V$ be such that $\phi'''(t)\geq 0$ for every $t\in [\min\{f(y) : y\sim x \text{ or } y=x\},f(x)]$.
Since $\phi$ is concave we have that $\phi(f(y))-\phi(f(x))\leq \phi'(f(x))(f(y)-f(x))$ for every $y\sim x$. If $f(y)<f(x)$ then we have moreover by Taylor's theorem with remainder that there exists $\xi_y\in (f(y),f(x))$ such that
\[
  \phi(f(y))-\phi(f(x)) = \phi'(f(x)) (f(y)-f(x)) + \frac{1}{2} \phi''(\xi_y) (f(y)-f(x))^2,
\]
and that $\phi''(\xi_y)\leq \phi'' (f(x))$ by the assumption on $\phi'''$. As such, we have the bound
\[
  \Delta(\phi\circ f)(x) \leq \phi'(f(x)) \Delta f(x)+ \frac{1}{2}\phi''(f(x)) \sum_{y\in V} P(x,y) (f(y)-f(x))^2 \mathbbm{1}(f(y)<f(x)),
\]
which is easily seen to imply the claim.
\end{proof}

We now apply \cref{lem:chain_rule} to prove \cref{lem:Laplacian_transformation_phi}. The proof is elementary calculus. Roughly speaking, the idea is to find solutions to the ODE
\[
  \lambda \phi'+\frac{1}{4d} \phi'' = \lambda'
\]
that are concave, increasing, and have $\phi'\leq 1$ and $\phi'''\geq 0$ on $[0,R]$, taking $\lambda'$ to be as large a negative number as possible with such solutions still existing. The relevant bounds on $\Delta (\phi \circ f)$ then following from \cref{lem:chain_rule} after extending $\phi$ to be constant on $[R,\infty)$. It will be convenient for the proof to slightly modify the above ODE, replacing the right hand side by a function less than $\lambda'$. (The proof will report on the outcome of this calculation rather than doing it from scratch.)

\begin{proof}[Proof of \cref{lem:Laplacian_transformation_phi}]
In each case below $\phi$ is constant on $[R,\infty)$, so that $\phi\circ f$ attains its maximum value at every $x$ with $f(x)\geq R$ and hence satisfies $\Delta(\phi\circ f)(x)\leq 0$ there; it therefore suffices to bound $\Delta(\phi\circ f)(x)$ for those $x\in X$ with $f(x)\leq R$.
If $\lambda\leq -1/4dR$ then
\[
  \lambda\leq \frac{-\lambda}{2\exp(4d\lambda R)-2} 
\]
and the claim holds with $\phi(t)=\min\{t,R\}$, since $\Delta(\phi\circ f)(x)\leq \Delta f(x)\leq \lambda$ whenever $f(x)\leq R$; so we may assume that $\lambda\geq -1/4dR$. First suppose that $\lambda> 0$. In this case, we take $\phi:[0,\infty)\to\R$ to be given by
\[
  \phi(t)= 
  \frac{1-e^{-4\lambda dt}-4\lambda dt e^{-4\lambda dR}}{4\lambda d (1-e^{-4\lambda dR})}
\]
for $t\in [0,R]$ and $\phi(t)=\phi(R)$ for $t>R$. One easily verifies that $\phi$ is concave and increasing with $\phi'''\geq 0$ and that $\phi$ solves the ODE
\[
  \lambda\phi'+\frac{1}{4d}\phi''=\frac{-\lambda}{e^{4\lambda dR}-1}
\]
on $[0,R]$. As such, it follows from \cref{lem:chain_rule} that
\[
  \Delta (\phi \circ f)(x) \leq \frac{-\lambda}{e^{4\lambda dR}-1} \mathbbm{1}(f(x)\leq R)
\]
for $x\in X$, which is stronger than the claim. Now suppose that $\lambda \in [-1/4dR,0]$. In this case, we take $\phi:[0,\infty)\to\R$ to be given by
\[
  \phi(t)=\frac{2Rt-t^2}{2R}
\]
for $t\in [0,R]$ and $\phi(t)=\phi(R)$ for $t>R$. As before, $\phi$ is easily verified to be concave and increasing with $\phi'''\geq 0$ (indeed $\phi'''\equiv 0$) and to solve the ODE
\[
  \lambda\phi'+\frac{1}{4d} \phi'' = \frac{\lambda(2R-2t)-1/2d}{2R} \leq -\frac{1}{4dR}
\]
on $[0,R]$, where the final inequality holds since $\lambda \in [-1/4dR,0]$.  As such, it follows from \cref{lem:chain_rule} that
\[
  \Delta (\phi \circ f)(x) \leq  -\frac{1}{4dR} \mathbbm{1}(f(x)\leq R),
\]
for all $x\in X$, which is stronger than the claim.
\end{proof}

\end{appendix}

\end{document}